\documentclass[twoside,notitlepage,12pt]{article}

\usepackage{amssymb}
\usepackage[leqno]{amsmath}
\usepackage{amsfonts}
\usepackage{amsopn}
\usepackage{amstext}
\usepackage{amsthm}

\usepackage[all]{xy}
\newdir{ >}{{}*!/-9pt/@{>}}

\usepackage{verbatim}
\usepackage[colorlinks, backref]{hyperref}

\newcommand{\En}{\mathcal N}
\newcommand{\Ha}{\mathcal H}

\providecommand{\cal}{\mathcal}
\renewcommand{\Bbb}{\mathbb}
\renewcommand{\frak}{\mathfrak}
\newenvironment{pf}{\begin{proof}}{\end{proof}}

\newcommand{\Aaa}{{\cal{A}}}
\newcommand{\Bee}{{\cal{B}}}

\newcommand{\Ef}{{\cal{F}}}
\newcommand{\Gee}{{\cal{G}}}

\newcommand{\Pee}{{\cal{P}}}

\newcommand{\Emm}{{\frak{M}}}

\newcommand{\Qyu}{{\Bbb{Q}}}
\newcommand{\Err}{{\Bbb{R}}}

\newcommand{\sig}{\sigma}
\newcommand{\eps}{\varepsilon}
\renewcommand{\phi}{\varphi}
\renewcommand{\rho}{\varrho}

\newcommand{\rest}{\restriction}

\newcommand{\ntr}{{n\in\omega}}

\newcommand{\loe}{\leqslant}
\newcommand{\goe}{\geqslant}

\newcommand{\subs}{\subseteq}
\newcommand{\sups}{\supseteq}
\newcommand{\nnempty}{\ne\emptyset}

\newcommand{\ovr}{\overline}
\newcommand{\til}{\tilde}
\renewcommand{\iff}{\Longleftrightarrow}

\newcommand{\id}[1]{{\operatorname{id}_{#1}}} 

\newcommand{\oraz}{\qquad\text{and}\qquad}

\newcommand{\Es}{{\cal{S}}}

\newcommand{\Ve}{{\mathbb V}}

\newtheorem{tw}{Theorem}[section]
\newtheorem{wn}[tw]{Corollary}
\newtheorem{lm}[tw]{Lemma}
\newtheorem{prop}[tw]{Proposition}

\theoremstyle{definition}
\newtheorem{df}[tw]{Definition}
\newtheorem{ex}[tw]{Example}

\theoremstyle{remark}

\newcommand{\setof}[2]{\{#1\colon #2\}}

\newcommand{\sett}[2]{\{#1\}_{#2}}
\newcommand{\sn}[1]{\{#1\}} 
\newcommand{\dn}[2]{\{#1,#2\}} 
\newcommand{\pair}[2]{\langle #1, #2 \rangle} 
\newcommand{\map}[3]{#1\colon #2 \to #3} 
\newcommand{\img}[2]{#1[#2]} 

\newcommand{\fra}{Fra\"iss\'e}
\newcommand{\jon}{J\'onsson}

\newcommand{\U}{\mathbb U}

\providecommand{\nat}{\omega}

\newcommand{\ciag}[1]{{\sett{{#1}_n}{\ntr}}}

\newcommand{\uzuple}[2]{{{\operatorname{Seq}}_{\loe{#1}}{\left(#2\right)}}}

\newcommand{\ciagi}{\sig}

\newcommand{\ciagipo}{\sig^{\operatorname{PO}}}

\newcommand{\norm}[1]{\|#1\|}

\newcommand{\clbal}{\overline{\operatorname{B}}}

\newcommand{\fK}{{\mathfrak{K}}}
\newcommand{\fL}{{\mathfrak{L}}}

\newcommand{\fR}{{\mathfrak{R}}}

\newcommand{\homos}[1]{ {}^{\operatorname{hom}}{#1} }
\newcommand{\embes}[1]{ {}^{\operatorname{emb}}{#1} }
\newcommand{\quos}[1]{ {}^{\operatorname{quo}}{#1} }

\newcommand{\los}[1]{{\mathfrak L\!\mathfrak O}_{ < {#1}}}
\newcommand{\flos}{\los\omega}

\newcommand{\bV}{{\mathbb{V}}}

\newcommand{\cmp}{\circ} 

\newcommand{\banach}{\ensuremath{\mathfrak B}} 
\newcommand{\banachi}{\ensuremath{\mathfrak B^{\operatorname{iso}}}} 

\newcommand{\metrics}{{\mathfrak M}} 

\newcommand{\haha}[1]{{\ensuremath{(\mathfrak H_{#1})}}}
\newcommand{\hahah}{{\ensuremath{(\mathfrak H)}}}

\newcommand{\separator}{\begin{center}***\end{center}}

\newcommand{\ob}[1]{\operatorname{Ob}\left(#1\right)}

\newcommand{\gur}{{\mathbb G}} 

\title{Injective objects and retracts of \fra\ limits}
\author{
{\sc Wies{\l}aw Kubi\'s}
\footnote{Research supported in part by the GA\v CR grant No. P 201/12/0290.}
\\ \\
{\small Mathematical Institute, Academy of Sciences of the Czech Republic}\\
and \\
{\small Institute of Mathematics,}
{\small Jan Kochanowski University in Kielce, Poland}\\
{\small\texttt{kubis@math.cas.cz}, \texttt{wkubis@pu.kielce.pl}}
}

\begin{document}

\maketitle

\begin{center}
\textit{Dedicated to the memory of my friend Pawe{\l} Waszkiewicz}
\end{center}

\begin{abstract}
We present a purely category-theoretic characterization of retracts of \fra\ limits. For this aim, we consider a natural version of injectivity with respect to a pair of categories (a category and its subcategory). It turns out that retracts of \fra\ limits are precisely the objects that are injective relatively to such a pair.
One of the applications is a characterization of non-expansive retracts of Urysohn's universal metric space.


\end{abstract}

\tableofcontents

\section{Introduction}

Throughout nearly all areas of mathematics one can find certain canonical objects that are uniquely determined by their homogeneity-like properties.
Historically, the first detected example of this sort was the set of rational numbers $\Qyu$, characterized by Cantor as the unique dense countable linearly ordered set with no end-points.
Another example is Urysohn's universal metric space, the unique separable complete metric space $\U$ containing isometric copies of all separable metric spaces and with the property that every isometry between finite subsets of $\U$ extends to a bijective isometry of $\U$.
In 1954, Roland \fra\ developed a general theory in the language of first-ordered structures, currently known as \emph{\fra\ theory}.
After his work, several universal homogeneous structures (called \emph{\fra\ limits}) had been identified and studied,
being important objects in various areas of mathematics, computer science and even mathematical physics~\cite{Droste}.
One needs to admit that the Urysohn space had been almost forgotten for many years, and not linked to \fra\ theory until a relatively recent line of research dealing with topological dynamics of automorphism groups.
A notable work in this area is~\cite{KPT}.
For more information on current status of \fra\ theory we refer to a recent survey article of Macpherson's~\cite{Macpherson}.

\separator

Recall that a \emph{\fra\ class} is a countable class $\Ef$ of finitely generated models of a fixed first-order language, satisfying the following conditions:
\begin{enumerate}
	\item[(i)] Given $a,b\in \Ef$ there exists $d\in \Ef$ such that both $a$ and $b$ embed into $d$ (Joint Embedding Property).
	\item[(ii)] Given $a,b\in \Ef$ and embeddings $\map ica$, $\map jcb$, there exist $w\in\Ef$ and embeddings $\map kaw$ and $\map \ell bw$ such that $k \cmp i = \ell \cmp j$ (Amalgamation Property).
	\item[(iii)] Given $a\in \Ef$, every substructure of $a$ is isomorphic to an element of $\Ef$.
\end{enumerate}
The \emph{\fra\ limit} of $\Ef$ is a countable model $U$ such that, up to isomorphism,
$$\Ef = \setof{a \subs U}{a \text{ is a finitely generated substructure of }U}$$
and for every isomorphism $\map hab$ between finitely generated substructures of $U$ there exists an automorphism $\map HUU$ such that $H \sups h$.
The latter property is called \emph{ultrahomogeneity}.
It is a classical theorem of \fra~\cite{F1} that the \fra\ limit exists and is unique, up to isomorphism.
Uncountable versions of \fra\ limits were studied by \jon~\cite{jon56, jon60}, supplemented by Morley and Vaught~\cite{MorleyVaught}.

A recent result of Dolinka~\cite{Dolinka} characterizes, under certain assumptions, countable models that are embeddable as retracts into the \fra\ limit.
Namely, he proves that, under certain conditions on the \fra\ class, retracts of the \fra\ limit are precisely the (countable) algebraically closed models.
Further study, in the context transformation semigroups and permutation group theory has been done in a recent PhD thesis of McPhee~\cite{McPhee}.

The aim of this note is to extend Dolinka's characterization to the case of category-theoretic \fra\ limits, at the same time weakening the assumption on the class of objects.
In particular, Dolinka's result assumes that models are finite and for each natural number $n$ there exist only finitely many isomorphic types of models generated by a set of cardinality $n$.
We do not make any of these assumptions.
Our result relates retracts of \fra\ limits to a natural variant of injectivity.
Among new applications, we characterize non-expansive retracts of the universal metric space of Urysohn. This metric space is formally not a \fra\ limit, because the category of finite metric spaces is uncountable. However, it can be ``approximated" by \fra\ limits of countable subcategories (e.g. by considering rational distances only).

Category-theoretic approach to \fra\ limits comes from the author's paper \cite{Kubfra}, motivated by a much earlier work of Droste and G\"obel \cite{DrGoe92} and by a recent work of Irwin and Solecki~\cite{IrSol} on projective \fra\ limits.
In~\cite{Kubfra} the key notion is a \emph{\fra\ sequence} rather than a \fra\ limit.
This turns out to be convenient, allowing to work in a single category (corresponding to finitely generated models), forgetting about the existence or non-existence of colimits.
In order to speak about retractions, we need to work with a pair of categories, both with the same objects; the first one allows ``embeddings" only, while the second one allows all possible homomorphisms.

\subsection{Categories of sequences}

Fix a category $\fK$. We shall treat $\fK$ as the class of arrows, the class of objects will be denoted by $\ob\fK$ and the set of $\fK$-arrows with domain $x$ and codomain $y$ will be denoted by $\fK(x,y)$.
A \emph{sequence} in $\fK$ is simply a covariant functor from $\omega$ into $\fK$.
One can think that the objects of $\fK$ are ``small" structures (e.g. finitely generated models of a fixed language).
Sequences in $\fK$ form a bigger category of ``large" structures.
For category-theoretic notions we refer to \cite{MacLane}.

We shall use the following convention: Sequences in $\fK$ will be denoted by capital letters $X,Y,Z,\dots$ and the objects of $\fK$ will be denoted by small letters $x,y,z,\dots$.
Fix a sequence $\map X\omega\fK$. 
Recall that formally $X$ assigns to each natural number $n$ an object $X(n)$ of $\fK$ and $X$ assigns a $\fK$-arrow $\map {X(n,m)}{X(n)}{X(m)}$ for each pair $\pair nm$ of natural numbers such that $n\loe m$.
We shall always write $x_n$ instead of $X(n)$ and $x_n^m$ instead of $X(n,m)$.
Note that being a functor imposes the conditions $x_n^n = \id{x_n}$ and $x_k^m = x_\ell^m \cmp x_k^\ell$ for $k\loe \ell \loe m$.

An arrow from a sequence $X$ to a sequence $Y$ is, by definition, a natural transformation from the functor $X$ into the functor $Y \cmp \psi$, where $\map \psi \omega\omega$ is increasing (i.e. $\psi$ is a covariant functor from $\omega$ to $\omega$).
We identify arrows that ``potentially converge" to the same limit.
More precisely, given natural transformations $\tau_0$ and $\tau_1$ from the sequence $X$ to $Y \cmp \psi_0$ and $Y \cmp \psi_1$, respectively, we say that $\tau_0$ is \emph{equivalent} to $\tau_1$, if the diagram consisting of both sequences $X$, $Y$ together with all arrows induced by $\tau_0$ and $\tau_1$ is commutative.
This is indeed an equivalence relation and it commutes with the composition, therefore $\ciagi\fK$ becomes a category.
In order to illustrate this idea, observe that every sequence is isomorphic to its cofinal subsequence.
Indeed, if $X$ is a sequence and $k = \ciag k$ is a strictly increasing sequence of natural numbers, then the $\ciagi\fK$-arrow $\map I{X\cmp k}X$ defined by $I=\ciag i$, where $i_n = \id{x_{k_n}}$, is an isomorphism.
Its inverse is $J = \ciag j$, where $j_n = x_n^{k_m}$ and $m=\min\setof{s}{k_s \goe n}$.
The composition $I \cmp J$ is formally $\ciag j$ regarded as an arrow from $X$ to $X$. Clearly, $I\cmp J$ is equivalent to the identity $\sett{\id{x_n}}{\ntr}$.
Similarly, $J\cmp I$ is equivalent to the identity of $X\cmp k$.

The original category $\fK$ may be regarded as a subcategory of $\ciagi\fK$, identifying an object $x$ with a sequence
$$\xymatrix{
x \ar[r]^{\id x} & x \ar[r]^{\id x} & x \ar[r]^{\id x} & \dots
}$$
Thus, we shall always assume that $\fK \subs \ciagi\fK$.
Given a sequence $X$ and $\ntr$, we shall denote by $x_n^\infty$ the arrow from $x_n$ to $X$ induced by the $n$th object of $X$.
Formally, $x_n^\infty$ is the equivalence class of $\sett{x_n^m}{m\goe n}$.

\subsection{\fra\ sequences}

\fra\ classes and limits can be described using categories. 
Let $\fK$ be a fixed category.

A \emph{\fra\ sequence} in $\fK$ is a sequence $U$ satisfying the following two conditions.
\begin{enumerate}
	\item[(F1)] For every object $x$ in $\fK$ there exist $\ntr$ and a $\fK$-arrow $x \to u_n$.
	\item[(F2)] For every $\ntr$ and for every $\fK$-arrow $\map f{u_n}y$ there exist $m \goe n$ and a $\fK$-arrow $\map gy{u_m}$ such that $g \cmp f = u_n^m$.
\end{enumerate}
Recall that $\fK$ has the \emph{amalgamation property} if for every $\fK$-arrows $\map fca$, $\map gcb$ there exist $\fK$-arrows $\map {f'}aw$, $\map {g'}bw$ satisfying $f' \cmp f = g' \cmp g$.
A \fra\ sequence exists whenever $\fK$ has the amalgamation property, the joint embedding property and has countably many isomorphic types of arrows.
A \fra\ sequence is unique up to isomorphism.
We refer to \cite{Kubfra} for the details.

A standard induction shows that the amalgamation property partially extends to the category of sequences. Namely:

\begin{prop}\label{prmknoteghb}
Assume $\fK$ has the amalgamation property. Then for every $\ciagi\fK$-arrows $\map fcA$, $\map gcB$ with $c\in \ob\fK$, there exist $\ciagi\fK$-arrows $\map {f'}AW$, $\map {g'}BW$ satisfying $f' \cmp f = g' \cmp g$.
\end{prop}

However, it is shown in~\cite{Kubfra} that in general the amalgamation property of $\fK$ does not imply the same property of $\ciagi\fK$.

Now, let $\fK\subs \fL$ be a pair of categories such that $\fK$ has the same objects as $\fL$. For instance, $\fK$ is a category of finitely generated models of a fixed language with embeddings and $\fL$ allows all homomorphisms.
Note that $\ciagi\fK$ as a subcategory of $\ciagi\fL$.
We shall need to deal with the category $\fR = \ciagi{(\fK,\fL)}$ whose objects are $\omega$-sequences in $\fK$ and the arrows come from $\fL$, i.e., $\ob\fR = \ob{\ciagi\fK}$ and $\fR(X,Y) = \ciagi\fL(X,Y)$ for $X,Y\in \ob{\fR}$.

For example, if $\fK$, $\fL$ are as above, $\ciagi{(\fK,\fL)}$ is the category of countable models with all possible homomorphisms, while $\ciagi\fK$ is the category of countable models with embeddings.

\section{Main result}

Let $\fK\subs \fL$ be two fixed categories with the same objects.
We say that $\pair \fK\fL$ has the \emph{mixed amalgamation property} if for every arrows $\map fca$ and $\map gcb$ such that $f\in \fK$ and $g\in \fL$, there exist arrows $\map {f'}aw$, $\map {g'}bw$ satisfying $f' \cmp f = g' \cmp g$ and such that $g' \in \fK$ and $f' \in \fL$.
The mixed amalgamation is described in the following diagram, where $\xymatrix{\ar@{ >->}[r] &}$ denotes an arrow in $\fK$.
$$\xymatrix{
a \ar[r]^{f'}& w\\
c \ar@{ >->}[u]^f \ar[r]_g & b \ar@{ >->}[u]_{g'}
}$$
We say that $\pair \fK\fL$ has the \emph{amalgamated extension property} if for every commutative $\fL$-diagram
$$\xymatrix{
a \ar[r]^f & x \\
c \ar@{ >->}[u]^i \ar@{ >->}[r]_j & b \ar[u]_g
}$$
with $i,j \in \fK$, there exist $\fK$-arrows $\map exy$, $\map kaw$, $\map \ell bw$ and an $\fL$-arrow $\map hwy$ such that $e \cmp f = h \cmp k$, $e \cmp g = h \cmp \ell$ and $k \cmp i  = \ell \cmp j$.
That is, the following diagram is commutative.
$$
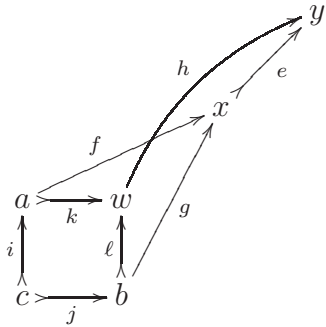
\xymatrix{
& & & y \\
& & x \ar@{ >->}[ur]_e & \\
a \ar[urr]^(.40)f \ar@{ >->}[r]_k & w \ar@/^1pc/[uurr]^{h} & & \\
c \ar@{ >->}[u]^i \ar@{ >->}[r]_j & b \ar[uur]_g \ar@{ >->}[u]^\ell & &
}$$

We now define the following axioms for a pair of categories $\pair \fK \fL$, needed for our main result.
\begin{enumerate}
	\item[\haha 0] $\fK \subs \fL$ and $\ob \fK = \ob \fL$.
	\item[\haha 1] $\fK$ has both the amalgamation property and the joint embedding property.
	\item[\haha 2] $\pair \fK \fL$ has the mixed amalgamation property.
	\item[\haha 3] $\pair \fK \fL$ has the amalgamated extension property.
\end{enumerate}

\begin{df}
A pair of categories $\pair \fK \fL$ \emph{has property \hahah}\ if it satisfies conditions \haha0 -- \haha3.
\end{df}

It is necessary to make some comments on the properties described above.
Namely, the condition $\ob\fK=\ob\fL$ can be removed from \haha0, it appears there for the sake of convenience only. The role of $\fL$ is offering more arrows than $\fK$, some of them will be needed for constructing retractions.
One can think of the $\fK$-arrows as ``embeddings". In most cases, these will be indeed monics.
Condition \haha1 is needed mainly for the existence and good properties of a \fra\ sequence in $\fK$.
Recall that the joint embedding property follows from amalgamations, whenever $\fK$ has an initial object (or at least a weakly initial object).
Condition \haha2 will be crucial for proving that the \fra\ sequence and its retracts are $\fK$-injective (see the definition below).
Finally, the somewhat technical condition \haha3 will be needed for the argument in the main lemma relating $\fK$-injective objects with the \fra\ sequence.
If $\fL$ has a terminal object then \haha3 implies that $\fK$ has the amalgamation property.
Summarizing, if $\fK$ has a weakly initial object and $\fL$ has a terminal object, then we may ignore condition \haha1.
Condition \haha3 becomes trivial if $\fK$ has pushouts in $\fL$.
We say that $\fK$ has \emph{pushouts in $\fL$} if for every pair of $\fK$-arrows $\map ica$, $\map jcb$, there exist $\fK$-arrows $\map kaw$, $\map \ell bw$ such that
$$\xymatrix{
a \ar@{ >->}[r]^k & w \\
c \ar@{ >->}[u]^i \ar@{ >->}[r]_j & b \ar@{ >->}[u]_\ell
}$$
is a pushout square in $\fL$.
It is obvious from the definition of a pushout that $\pair\fK\fL$ has the amalgamated extension property (with $y=x$ and $e = \id x$) whenever $\fK$ has pushouts in $\fL$.
Let us remark that for all examples with property \hahah\ appearing in this note, the amalgamated extension property holds with $x=y$ and $e=\id x$ (see the definition and diagram above).

Below is the crucial notion, whose variations appear often in the literature (see, e.g., \cite{AHRT}, where a definition similar to ours can be found).

\begin{df}
Let $\fK \subs \fL$ be two categories with the same objects.
We say that $A\in \ob {\ciagi\fK}$ is \emph{$\fK$-injective in} $\ciagi{(\fK,\fL)}$ if for every $\fK$-arrow $\map i a b$, for every $\ciagi{(\fK,\fL)}$-arrow $\map f a A$, there exists a $\ciagi{(\fK,\fL)}$-arrow $\map {\ovr f} b A$ such that ${\ovr f} \cmp i = f$.
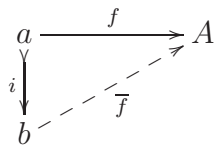
$$\xymatrix{
a \ar@{ >->}[d]_i \ar[rr]^f & & A \\
b \ar@{-->}[rru]_{\ovr f}
}$$
\end{df}

This definition obviously generalizes to an arbitrary pair of categories $\fK \subs \fR$. We restrict attention to the special case $\fR = \ciagi{(\fK,\fL)}$, since more general versions will not be needed.

Following is a useful criterion for injectivity.

\begin{prop}\label{pinkarerg}
Assume $\pair \fK\fL$ has the mixed amalgamation property and $X \in \ob{\ciagi\fK}$. Then $X$ is $\fK$-injective in $\ciagi{(\fK,\fL)}$ if and only if
for every $\ntr$, for every $\fK$-arrow $\map f{x_n}y$, there exist $m\goe n$ and an $\fL$-arrow $\map gy{x_m}$ satisfying
$$g \cmp f = x_n^m.$$
\end{prop}

\begin{pf}
Suppose $X$ is $\fK$-injective and fix a $\fK$-arrow $\map f{x_n}y$.
Applying $\fK$-injectivity for $\map {x_n^\infty}{x_n}X$, we find $\map GyX$ such that $G \cmp f = x_n^\infty$.
The arrow $G$ factors through some $\fL$-arrow $\map gy{x_m}$ for some $m\goe n$, that is, $G = x_m^\infty \cmp g$.
Finally, $g \cmp f = x_n^m$.

Suppose now that $X$ satisfies the condition above and fix a $\fK$-arrow $\map jab$ and a $\ciagi{(\fK,\fL)}$-arrow $\map Fa{X}$.
Then $F = x_n^\infty \cmp f$ for some $\fL$-arrow $f$, where $\ntr$.
Applying the mixed amalgamation property, find a $\fK$-arrow $\map h{x_n}y$ and an $\fL$-arrow $\map gby$ such that $g \cmp j = h \cmp f$.
By assumption, there exist $m \goe n$ and an $\fL$-arrow $\map ky{x_m}$ such that the following diagram commutes.
$$
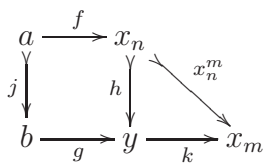
\xymatrix{
a \ar@{ >->}[d]_j \ar[r]^f & x_n \ar@{ >->}[d]_h \ar@{ >->}[dr]^{x_n^m} & \\
b \ar[r]_g & y \ar[r]_k & x_m
}$$
Finally, taking $G = x_m^\infty \cmp k \cmp g$, we get $G \cmp j = F$.
\end{pf}

Our interest in $\fK$-injectivity comes from the following fact, which is an immediate consequence of the criterion above.

\begin{prop}\label{pwkindrzektif}
Assume $\pair \fK\fL$ has the mixed amalgamation property and $U$ is a \fra\ sequence in $\fK$. Then $U$ is $\fK$-injective in $\ciagi{(\fK,\fL)}$.
\end{prop}

We shall need the following ``injective" version of amalgamated extension property.

\begin{lm}\label{lekstaminja}
Assume $\pair \fK\fL$ satisfies \hahah\ and $X\in \ob{\ciagi\fK}$ is $\fK$-injective in $\ciagi{(\fK,\fL)}$.
Then for every $\fK$-arrows $\map ica$, $\map jcb$ and for every $\ciagi{(\fK,\fL)}$-arrows $\map FaX$, $\map GbX$ such that $F \cmp i = G \cmp j$, there exist $\fK$-arrows $\map kaw$, $\map \ell bw$ and a $\ciagi{(\fK,\fL)}$-arrow $\map HwX$ such that the diagram
$$
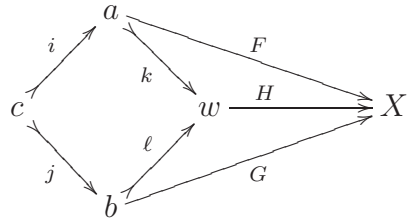
\xymatrix{
& a \ar@{ >->}[rd]_k \ar[rrrd]^F & & & \\
c \ar@{ >->}[ru]^i \ar@{ >->}[rd]_j &  & w \ar[rr]^(.30)H & & X\\
& b \ar@{ >->}[ru]^\ell \ar[rrru]_G & & &
}$$
commutes.
\end{lm}

\begin{pf}
Find $n$ such that $F = x_n^\infty \cmp f$ and $G = x_n^\infty \cmp g$ for some $\fL$-arrows $f,g$, where $\map {x_n^\infty}{x_n}X$ is the canonical arrow induced by the $n$th object of the sequence $X$.
Using property \haha3, we find $\fK$-arrows $\map kaw$, $\map \ell bw$, $\map e{x_n}y$ and an $\fL$-arrow $\map hwy$ such that $h \cmp k = e \cmp f$ and $h \cmp \ell = e \cmp g$.
Using the $\fK$-injectivity of $X$ we can find a $\ciagi{(\fK,\fL)}$-arrow $\map PyX$ such that $P \cmp e = x_n^\infty$.
Let $H = P \cmp h$.
Then
$$H \cmp k = P \cmp h \cmp k = P \cmp e \cmp f = x_n^\infty \cmp f = F.$$
Similarly, $H \cmp \ell = G$.
\end{pf}

The following lemma is crucial.

\begin{lm}\label{lkrusall}
Assume $\pair \fK\fL$ has property \hahah\ and $A$ is a $\fK$-injective object in $\ciagi{(\fK,\fL)}$.
Furthermore, assume $U$ is a \fra\ sequence in $\fK$ and $\map FXA$ is an arbitrary $\ciagi{(\fK,\fL)}$-arrow.
Then there exist a $\ciagi\fK$-arrow $\map JXU$ and a $\ciagi{(\fK,\fL)}$-arrow $\map GUA$ such that $G \cmp J = F$.
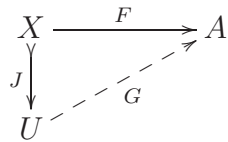
$$\xymatrix{
X \ar@{ >->}[d]_J \ar[rr]^F & & A \\
U \ar@{-->}[rru]_{G}
}$$
\end{lm}

\begin{pf}
Recall that we use the usual convention for objects $x_n = X(n)$, $u_n = U(n)$, and for arrows $x^m_n = X(n,m)$, $u^m_n = U(n,m)$.
We shall construct inductively the following ``triangular matrix" in $\fK$, together with commuting $\ciagi{(\fK,\fL)}$-arrows $\map {F_{i,j}}{w_{i,j}}A$ for $j\loe i+1$, where we agree that $w_{i,0} = x_i$ and $w_{i,i+1} = u_{\ell_i}$.
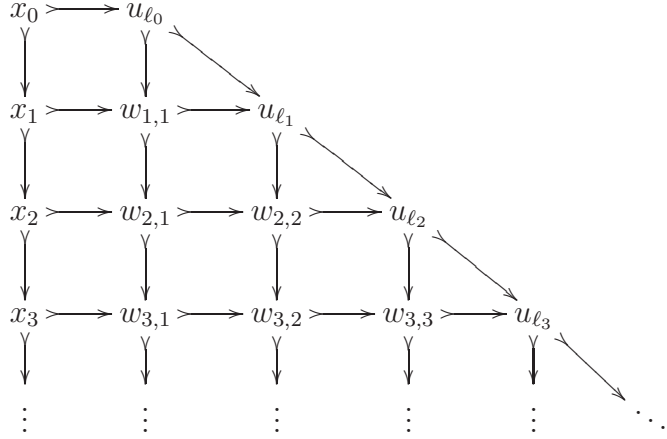
$$\xymatrix{
x_0\ar@{ >->}[d]\ar@{ >->}[r] & u_{\ell_0}\ar@{ >->}[d]\ar@{ >->}[dr] & & & & \\
x_1\ar@{ >->}[d]\ar@{ >->}[r] & w_{1,1}\ar@{ >->}[d]\ar@{ >->}[r] & u_{\ell_1}\ar@{ >->}[d]\ar@{ >->}[dr] & & & \\
x_2\ar@{ >->}[d]\ar@{ >->}[r] & w_{2,1}\ar@{ >->}[d]\ar@{ >->}[r] & w_{2,2}\ar@{ >->}[d]\ar@{ >->}[r] & u_{\ell_2}\ar@{ >->}[d]\ar@{ >->}[dr] & & \\
x_3\ar@{ >->}[d]\ar@{ >->}[r] & w_{3,1}\ar@{ >->}[d]\ar@{ >->}[r] & w_{3,2}\ar@{ >->}[d]\ar@{ >->}[r] & w_{3,3}\ar@{ >->}[d]\ar@{ >->}[r] & u_{\ell_3}\ar@{ >->}[d]\ar@{ >->}[dr] & \\
\vdots & \vdots & \vdots & \vdots & \vdots & \ddots \\
}$$
The first column in the diagram above is the sequence $X$, while the diagonal is a cofinal subsequence of $U$.
Our initial assumption on $F_{i,j}$ is that $\sett{F_{n,0}}{\ntr} = F$.
It is clear how to start the construction: Using the \fra\ property of $U$, we find $\ell_0$ and a $\fK$-arrow $\map {e_0}{x_0}{u_{\ell_0}}$.
Next, using the $\fK$-injectivity of $A$, we find $\map{F_{0,1}}{u_{\ell_0}}A$ satisfying $F_{0,1} \cmp e_0 = F_{0,0}$.

Suppose the $n$th row has already been constructed, together with arrows $F_{i,j}$ for $i\loe n$, $j\loe n+1$.
Starting from $\fK$-arrows $\map {x_n^{n+1}}{x_n}{x_{n+1}}$ and $x_n \to w_{n,1}$, using Lemma~\ref{lekstaminja},
we find $w_{n+1,1}\in \ob\fK$ and $\fK$-arrows $w_{n,1} \to w_{n+1,1}$, $x_{n+1} \to w_{n+1,1}$, and a $\ciagi{(\fK,\fL)}$-arrow $\map {F_{n+1,1}}{w_{n+1,1}}A$ such that the diagram
$$
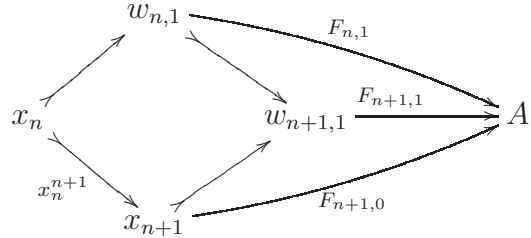
\xymatrix{
& w_{n,1}\ar@{ >->}[rd] \ar@/^/[rrrd]^{F_{n,1}} & & & \\
x_n\ar@{ >->}[ru]\ar@{ >->}[rd]_{x_n^{n+1}} & & w_{n+1,1} \ar[rr]^(.4){F_{n+1,1}} & & A\\
& x_{n+1}\ar@{ >->}[ru] \ar@/_/[rrru]_{F_{n+1,0}} & & &
}$$
commutes.
Continuing this way, using Lemma~\ref{lekstaminja}, we obtain the $(n+1)$st row and $\ciagi{(\fK,\fL)}$-arrows $F_{n+1,i}$ for $i\loe n+1$ which commute together with the following diagram.
$$
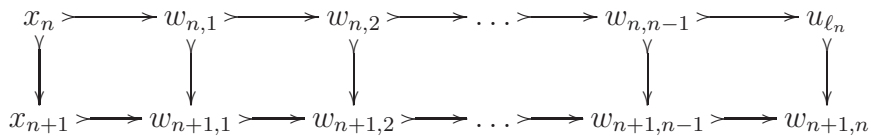
\xymatrix{
x_n\ar@{ >->}[d]\ar@{ >->}[r] & w_{n,1}\ar@{ >->}[d]\ar@{ >->}[r] & w_{n,2}\ar@{ >->}[d]\ar@{ >->}[r] & \dots\ar@{ >->}[r] & w_{n,n-1}\ar@{ >->}[d]\ar@{ >->}[r] & u_{\ell_n}\ar@{ >->}[d] \\
x_{n+1}\ar@{ >->}[r] & w_{n+1,1}\ar@{ >->}[r] & w_{n+1,2}\ar@{ >->}[r] & \dots\ar@{ >->}[r] & w_{n+1,n-1}\ar@{ >->}[r] & w_{n+1,n}
}$$
Now, using the \fra\ property of $U$ we find $\ell_{n+1} > \ell_n$ and a $\fK$-arrow $w_{n+1,n} \to u_{\ell_{n+1}}$ making the triangle
$$\xymatrix{
u_{\ell_n}\ar@{ >->}[d]\ar@{ >->}[dr] &  \\
w_{n+1,n}\ar@{ >->}[r] & u_{\ell_{n+1}}
}$$
commutative.
Using Lemma~\ref{lekstaminja} again, we get an arrow $\map{F_{n+1,n+2}}{u_{\ell_{n+1}}}A$ commuting with $F_{n+1,n}$, $F_{n,n+1}$ and the triangle above.

Finally, the compositions of the horizontal arrows in the triangular ``matrix" constructed above induce an arrow of sequences $\map JXU$ in $\ciagi\fK$.
The inductive construction also gives a sequence of arrows $\sett{F_{n,n+1}}{\ntr}$ that turns into a $\ciagi{(\fK,\fL)}$-arrow $\map GUA$ satisfying $G \cmp J = F$.
This completes the proof.
\end{pf}

\begin{tw}\label{tmejnn}
Let $\pair \fK \fL$ be a pair of categories with property \hahah. 
Assume $\fK$ has a \fra\ sequence $U$ and let $X$ be an arbitrary sequence in $\fK$.
The following properties are equivalent.
\begin{enumerate}
	\item[(a)] $X$ is $\fK$-injective in $\ciagi{(\fK, \fL)}$.
	\item[(b)] There exist a $\ciagi{\fK}$-arrow $\map JXU$ and a $\ciagi{(\fK,\fL)}$ arrow $\map RUX$ such that $R \cmp J = \id X$.
	\item[(c)] $X$ is a retract of $U$ in $\ciagi{(\fK, \fL)}$.
\end{enumerate}
\end{tw}

Note that condition (c) is formally weaker than (b), since it is not required in (c) that the right inverse of a retraction $\map RUX$ is a $\ciagi\fK$-arrow.

\begin{pf}
(a) $\implies$ (b)
Applying Lemma~\ref{lkrusall} to the identity $\map {\id X}XX$, we get a $\ciagi\fK$-arrow $\map JXU$ and a $\ciagi{(\fK,\fL)}$-arrow $\map RUX$ such that $R \cmp J = \id X$.

(b) $\implies$ (c) This is obvious.

(c) $\implies$ (a)
Let $\map JXU$ and $\map RUX$ be $\ciagi{(\fK,\fL)}$-arrows such that $R \cmp J = \id X$.

Fix a $\fK$-arrow $\map iab$ and an $\fL$-arrow $\map FaX$.
By Proposition~\ref{pwkindrzektif}, $U$ is $\fK$-injective in $\ciagi{(\fK,\fL)}$, so there exists $\map GbU$ such that $G \cmp i = J \cmp F$.
Finally, we have $R \cmp G \cmp i = R \cmp J \cmp F = F$.
\end{pf}

\subsection{Remarks on absolute retracts}

One can have a false impression after reading our result characterizing retracts of a \fra\ sequence, namely that \emph{every} embedding (i.e. a $\ciagi\fK$-arrow) of a $\fK$-injective $\ciagi \fK$-object into the \fra\ sequence admits a left inverse in $\ciagi(\fK,\fL)$.
This is not true in general.
We make a brief discussion of this problem.
To be more concrete, we assume $\pair \fK \fL$ has property \hahah\ and let $U$ be a \fra\ sequence in $\fK$.
The problem stated above is strictly related to the following well known concept:

\begin{df}
Let $\fK \subs \fL$ be as above. We say that $W \in \ob{\ciagi \fK}$ is an \emph{absolute retract} in $\ciagi(\fK,\fL)$ if for every $\ciagi \fK$-arrow $\map J W Y$ there exists a $\ciagi(\fK,\fL)$-arrow $\map R Y W$ such that $R \cmp J = \id W$.
\end{df}

This notion is well known especially in topology.
In particular, there is a rich theory of absolute retracts in geometric topology (see \cite{Borsuk}).
One of the main aspects is the existence of some ``canonical" objects which can be used for checking whether a given object is an absolute retract or not.
For instance, in the category of compact topological spaces, an absolute retract is simply a retract of a Tikhonov cube.
In the category of metric spaces with continuous maps, absolute retracts are retracts of convex sets in normed linear spaces.
In the category of metric spaces with non-expansive maps, absolute retracts are hyperconvex metric spaces \cite{AroPan}.

Notice that our definition is relative to a fixed category of ``small" objects, e.g. spaces of weight less than a fixed cardinal number.
In the case of compact topological (or metric) spaces, the ``canonical" objects (e.g. Tikhonov cubes) turn out to be absolute retracts without restrictions on the weight of the spaces and therefore being an absolute retract in a category of objects of restricted ``size" is equivalent to being an absolute retract in the big category, with no ``size" restrictions.

The following fact is rather standard.

\begin{prop}\label{Pierfbrw}
Assume $\fK \subs \fL$ is a pair of categories such that $\pair{\ciagi \fK}{\ciagi(\fK,\fL)}$ has the mixed amalgamation property.
Given $W \in \ob{\ciagi\fK}$ the following two properties are equivalent.
\begin{enumerate}
	\item[(a)] $W$ is an absolute retract in $\ciagi(\fK,\fL)$.
	\item[(b)] $W$ is $\ciagi\fK$-injective in $\ciagi(\fK,\fL)$.
\end{enumerate}
\end{prop}

\begin{pf}
Only (a)$\implies$(b) requires an argument.
Fix a $\ciagi\fK$-arrow $\map I X Y$ and a $\ciagi(\fK,\fL)$-arrow $\map F X W$.
Using the mixed amalgamation property we find a $\ciagi\fK$-arrow $\map J W V$ and a $\ciagi(\fK,\fL)$-arrow $\map G Y W$ for which the diagram
$$\xymatrix{
W \ar@{ >->}[r]^J & V \\
X \ar[u]^F \ar@{ >->}[r]_I & Y \ar[u]_G
}$$
is commutative.
Let $\map R V W$ be such that $R \cmp J = \id W$.
Then $R \cmp G \cmp I = F$.
\end{pf}

Now the problem arises whether $\fK$-injectivity implies $\ciagi\fK$-injectivity.
The next result characterizes this property using a \fra\ sequence.

\begin{tw}
Let $\fK \subs \fL$ be a pair of categories with property \hahah, and let $U\in \ob {\ciagi \fK}$ be a \fra\ sequence in $\fK$.
Assume further that $\pair{\ciagi \fK}{\ciagi(\fK,\fL)}$ has the mixed amalgamation property.
Then the following statements are equivalent.
\begin{enumerate}
	\item[(a)] $\fK$-injectivity implies $\ciagi\fK$-injectivity in $\ciagi(\fK,\fL)$.
	\item[(b)] $U$ is $\ciagi\fK$-injective in $\ciagi(\fK,\fL)$.
	\item[(c)] For every $\ciagi\fK$-arrow $\map J U U$ there exists a $\ciagi(\fK,\fL)$-arrow $\map R U U$ such that $R \cmp J = \id U$.
\end{enumerate}
\end{tw}

\begin{pf}
Implication (a)$\implies$(b) follows from the fact that $U$ is $\fK$-injective (Proposition~\ref{pwkindrzektif}).
Implication (b)$\implies$(c) is trivial.
Implication (b)$\implies$(a) follows directly from Theorem~\ref{tmejnn}, because a retract of a $\ciagi\fK$-injective object is obviously $\ciagi\fK$-injective. It remains to show that (c)$\implies$(b).

Suppose $U$ is not $\ciagi\fK$-injective.
By Proposition~\ref{Pierfbrw}, there exists a $\ciagi\fK$-arrow $\map I U Y$ which is not left-invertible in $\ciagi(\fK,\fL)$. Here we have used the mixed amalgamation property for $\pair {\ciagi\fK}{\ciagi(\fK,\fL)}$.
Since $U$ is \fra, there exists a $\ciagi\fK$-arrow $\map J Y U$.
Using (c), we find a $\ciagi(\fK,\fL)$-arrow $\map R U U$ such that $R \cmp J \cmp I = \id U$.
But now $R \cmp J$ is a left inverse to $I$, a contradiction.
\end{pf}

An interesting consequence of the result above is that whenever $\fK$-injectivity is different from $\ciagi\fK$-injectivity, it is witnessed by some $\ciagi\fK$-arrow $\map J U U$ with no left inverse in $\ciagi(\fK,\fL)$.
In other words, $U$ carries all the information about $\ciagi\fK$-injectivity.

\subsection{Extensions of the main result}

Theorem~\ref{tmejnn} has a natural generalization to uncountable \fra\ sequences.
More precisely, let $\kappa$ be an uncountable regular cardinal and assume that all sequences in $\fK$ of length $ < \kappa$ have colimits in $\fL$, where the colimiting cocones are $\fK$-arrows.
In this case we say that $\fK$ is \emph{$\kappa$-continuous} in $\fL$.
Under this assumption, a version of Lemma~\ref{lkrusall} for $\kappa$-sequences is true, with almost the same proof---usual induction is replaced by transfinite induction.
Proposition~\ref{pwkindrzektif} is valid for arbitrary \fra\ sequences, the countable length of the sequence was never used in the proof.

Let $\uzuple \kappa\fK$ denote the category of all sequences in $\fK$ of length $\loe \kappa$, with arrows induced by natural transformations (like in the countable case).
Let $\uzuple \kappa{\fK,\fL}$ denote the category with the same objects as $\uzuple \kappa\fK$, and with arrows taken from $\uzuple \kappa\fL$.
We can now formulate an ``uncountable" version of our main result.

\begin{tw}\label{tmejnnunct}
Let $\kappa$ be an uncountable regular cardinal and let $\pair \fK \fL$ be a pair of categories with property \hahah, such that $\fK$ is $\kappa$-continuous in $\fL$. 
Assume $\fK$ has a \fra\ sequence $U$ of length $\kappa$.
Given a sequence $X$ in $\fK$ of length $\loe \kappa$, the following properties are equivalent.
\begin{enumerate}
	\item[(a)] $X$ is $\fK$-injective in $\uzuple \kappa{\fK, \fL}$.
	\item[(b)] There exist a $\uzuple \kappa{\fK}$-arrow $\map JXU$ and a $\uzuple \kappa{\fK,\fL}$ arrow $\map RUX$ such that $R \cmp J = \id X$.
	\item[(c)] $X$ is a retract of $U$ in $\uzuple \kappa{\fK, \fL}$.
\end{enumerate}
\end{tw}

Let us now come back to the countable case.
Assume $\pair \fK\fL$ has property \hahah\ and moreover $\fK$ has pushouts in $\fL$.
Let us look at the proof of Lemma~\ref{lkrusall}.
We can assume that all squares in the infinite ``triangular matrix" constructed there are pushouts in $\fL$.
Using the notation from the proof of Lemma~\ref{lkrusall}, let $W^n$ denote the sequence coming from the $n$th column.
Observe that the arrow from $W^n$ to $W^{n+1}$ is determined by the ``horizontal" $\fK$-arrow $w_{n+1,n} \to w_{n+1,n+1}$.
In other words, all other $\fK$-arrows come as a result of the corresponding pushout square.
An arrow of sequences $\map FVW$ determined by pushouts from a single $\fK$-arrow will be called \emph{pushout generated} from $\fK$.
Denote by $\ciagipo \fK$ the category whose objects are $\omega$-sequences in $\fK$, while arrows are pushout generated from $\fK$.
A deeper analysis of the proof of Lemma~\ref{lkrusall} gives the following observation, which may be of independent interest.

\begin{prop}
Assume $\pair \fK \fL$ is a pair of categories with property \hahah\ and $\fK$ has pushouts in $\fL$.
Let $X \in \ob{\ciagi\fK}$ be $\fK$-injective in $\ciagi{(\fK,\fL)}$. Then:
\begin{enumerate}
	\item[(1)] $X$ is $\ciagipo \fK$-injective in $\ciagi{(\fK,\fL)}$.
	\item[(2)] Let $U$ be a \fra\ sequence in $\fK$. There exists a sequence
$$X_0 \to X_1 \to X_2 \to \dots$$
in $\ciagipo \fK$ such that $X_0 = X$ and $U$ is the colimit of this sequence in $\ciagi\fK$.
\end{enumerate}
\end{prop}

Clearly, (1) and (2) imply immediately that $X$ is a retract of $U$.

\section{Applications}

We start with some more comments on property \hahah.
In many cases (especially in model-theoretic categories), it is much easier to prove the (mixed) amalgamation property for special ``primitive" arrows rather than for arbitrary arrows.
In order to formalize this idea, fix a pair of categories $\pair\fK\fL$ satisfying condition \haha0 and fix a collection $\Ef \subs \fK$ (actually $\Ef$ might be a proper class). We say that $\fK$ is \emph{generated} by $\Ef$ if for every $f\in \fK$ there exist $\ntr$ and $g_0,\dots, g_{n-1} \in \Ef$ such that
$f = g_{n-1} \cmp \dots \cmp g_0$.
For example, if $\fK$ is the category of embeddings of finite models of a fixed first-order language, $\Ef$ may be the class of embeddings $\map fST$ such that $T$ is generated by $\img fS \cup \sn b$ for some $b\in T$.
We define the amalgamation property for $\Ef$ and the mixed amalgamation property for $\pair\Ef\fL$, as before.

\begin{prop}\label{pprimitivv}
Let $\fK \subs \fL$ be two categories with the same objects, where $\fK$ has the joint embedding property.
Assume further that $\fK$ is generated by a family $\Ef$ such that $\Ef$ has the amalgamation property and $\pair \Ef\fL$ has both the mixed amalgamation property and the amalgamated extension property.
Then $\pair\fK\fL$ has property \hahah.
\end{prop}

\begin{pf}
Given an arrow $f\in \fK$, we say that $f$ \emph{has length} $\loe n$ if $f = g_{n-1} \cmp \dots \cmp g_0$, where $g_0, \dots, g_{n-1} \in \Ef$.
In particular, all arrows in $\Ef$ have length $1$.
Easy induction shows that if $\map ica$, $\map jcb$ are $\fK$-arrows such that the length of $i$ is $\loe m$ and the length of $j$ is $\loe n$, then there exist $\fK$-arrows $\map kaw$, $\map \ell bw$ such that $k \cmp i = \ell \cmp j$ and $k$ has length $\loe n$, while $\ell$ has length $\loe m$.
Since every $\fK$-arrow has a finite length, this shows that $\fK$ has the amalgamation property.

A similar induction on the length of $\fK$-arrows shows that $\pair \fK\fL$ has the amalgamated extension property.
Finally, using the fact that $\pair\Ef\fL$ has the mixed amalgamation property, we prove by induction that for every $\fK$-arrow $\map ica$ of length $\loe n$, and for every $\fL$-arrow $\map fcb$, there exist an $\fL$-arrow $\map gaw$ and a $\fK$-arrow $\map \ell bw$ of length $\loe n$ such that $g \cmp i = \ell \cmp f$. This shows that $\pair \fK\fL$ has the mixed amalgamation property.
\end{pf}

Another simplification for proving property \hahah\ is the concept of mixed pushouts.

Let $\fK\subs \fL$ be two categories with the same objects. We say that $\pair\fK\fL$ has the \emph{mixed pushout property} if for every arrows $\map fca$ and $\map gcb$ such that $f\in \fK$ and $g\in \fL$, there exist arrows $\map {f'}aw$ and $\map {g'}bw$ such that $f'\in \fL$, $g'\in \fK$ and 
$$\xymatrix{
a \ar[r]^{f'} & w \\
c \ar@{ >->}[u]^f \ar[r]_g & b \ar@{ >->}[u]_{g'}
}$$
is a pushout square in $\fL$.
Note that if both $f,g$ are $\fK$-arrows in the definition above, then so are $f',g'$, by uniqueness of the pushout.

The definition above makes sense (and is applicable) in case where $\fK$ is an arbitrary family of arrows, not necessarily a subcategory. This is presented in the next statement.

\begin{prop}\label{pjetghsd}
Let $\fK\subs \fL$ be two categories with the same objects. 
Assume that $\fK$ has the joint embedding property and $\Ef\subs \fK$ is such that $\pair \Ef\fL$ has the mixed pushout property and $\Ef$ generates $\fK$.
Then $\pair \fK\fL$ has property \hahah.
\end{prop}

\begin{pf}
Suppose first that $\Ef = \fK$.
The amalgamation property (condition \haha1) follows from the remark above, namely that the pushout of two $\fK$-arrows consists of $\fK$-arrows.
Mixed amalgamation property (condition \haha2) is just a weaker version of the mixed pushout property.
Finally, the amalgamated extension property (condition \haha3) follows immediately from the definition of a pushout.

Suppose now that $\Ef \ne \fK$. It suffices to prove that $\pair\fK\fL$ has the mixed pushout property.
Like in the proof of Proposition~\ref{pprimitivv}, we use induction on the length of $\fK$-arrows, bearing in mind that the obvious composition of two pushout squares is a pushout square.
More precisely, the inductive hypothesis says: Given a $\fK$-arrow $\map ica$ of length $< n$, and an $\fL$-arrow $\map fcb$, there exist an $\fL$-arrow $\map gaw$ and a $\fK$-arrow $\map \ell bw$ of length $< n$ such that
$$\xymatrix{
b \ar@{ >->}[r]^\ell & w \\
c \ar[u]^f \ar@{ >->}[r]_i & a \ar[u]_g
}$$
is a pushout square in $\fL$.
\end{pf}

Many natural pairs of categories, in particular coming from model theory, have the mixed pushout property. Concrete well known examples are finite graphs, partially ordered sets, semilattices. Each of these classes is considered as a pair of two categories, the first one with embeddings and the second one with all homomorphisms.
These examples are mentioned in \cite{Dolinka}.
A typical example of a pair $\pair \fK\fL$ with property \hahah, failing the mixed pushout property is the category $\fL$ of all finite linear orders with increasing (i.e. order preserving) functions and $\fK$ the category of all finite linear orders with embeddings.

In contrast to the above results, it is worth mentioning a \fra\ class that does not fit into our framework.
Namely, the \fra\ class of finite $K_n$-free graphs (where $K_n$ denotes the complete graph with $n$ vertices and $n > 2$) has the pushout property (formally the class of embeddings has pushouts in the class of all homomorphisms), yet the corresponding pair of categories fails to have mixed amalgamations.
Specifically, a graph is meant to be a structure with one symmetric irreflexive binary relation, so a homomorphism of graphs cannot identify vertices connected by edges.
In other words, every graph homomorphism restricted to a complete subgraph becomes an embedding.
It has been proved by Mudrinski~\cite{Mudrinski} that for $n > 2$, the \fra\ limit of $K_n$-free graphs (called the \emph{Henson graph} $H_n$) is retract rigid, i.e. identity is the only retraction of $H_n$.
On the other hand, we have the following easy fact (stated in a different form in \cite[Example 3.3]{DolinkaBerg}).

\begin{prop}
No $K_n$-free graph with $n > 2$ is injective for finite $K_n$-free graphs.
\end{prop}

\begin{pf}
Suppose $X$ is such a graph.
Using injectivity for $S=\emptyset$ and $T = K_{n-1}$, we see that $X$ contains an isomorphic copy $K$ of $K_{n-1}$.
Now let $S$ be a graph with $n-1$ vertices and no edges and let $\map fSX$ be a bijection onto $K$.
Let $T = S \cup \sn v$, where $v$ is connected to all the vertices of $S$.
By injectivity, there exists a homomorphism $\map gTX$ extending $f$.
But now $K \cup \sn{g(v)} \subs X$ is a copy of $K_n$, a contradiction.
\end{pf}

Before discussing concrete examples of pairs with property \hahah, we make one more remark on injectivity. Recall that an arrow $\map jxy$ is \emph{left-invertible} in $\fL$ if there exists $f\in \fL$ such that $f \cmp j = \id x$.
The following is an easy consequence of our main result.

\begin{wn}\label{wnrhojht}
Let $\pair \fK\fL$ be a pair of categories such that every $\fK$-arrow is left-invertible in $\fL$.
Assume that $\pair \fK\fL$ has property \hahah\ and $U$ is a \fra\ sequence in $\fK$.
Then for every sequence $X \in \ob {\ciagi\fK}$ there exist a $\ciagi\fK$-arrow $\map JXU$ and a $\ciagi{(\fK,\fL)}$-arrow $\map RUX$ such that $R \cmp J = \id X$.
\end{wn}

\begin{pf}
In view of Theorem~\ref{tmejnn}, it suffices to show that every sequence is $\fK$-injective in $\ciagi{(\fK,\fL)}$.
Fix $X \in \ob {\ciagi\fK}$, a $\fK$-arrow $\map jab$, and a $\ciagi{(\fK,\fL)}$-arrow $\map faX$.
Choose an $\fL$-arrow $\map rba$ such that $r \cmp j = \id a$.
Then $g = f \cmp r$ has the property that $g \cmp j = f$.
This shows that $X$ is $\fK$-injective in $\ciagi{(\fK,\fL)}$.
\end{pf}

This corollary applies to finite Boolean algebras (also noted in \cite{Dolinka}) and, as we shall see later, to finite linear orderings.

\subsection{\fra\ classes and algebraically closed models}

Let $\Emm$ be a class of finitely generated models of a fixed first-order language $L$.
It is natural to consider the category $\homos\Emm$ whose objects are all elements of $\Emm$ and arrows are all homomorphisms (i.e. maps that preserve all relations, functions and constants).
It is also natural to consider the category $\embes\Emm$ whose objects are again all elements of $\Emm$, while arrows are embeddings only.
In many cases, $\pair{\embes\Emm}{\homos\Emm}$ has property \hahah.

Simplifying the notation, we shall say that $\Emm$ has the {\em pushout property} or {\em mixed amalgamation property} if $\pair{\embes\Emm}{\homos\Emm}$ has such a property.
Denote by $\ovr \Emm$ the class of all (countable) models that are unions of $\omega$-chains of models from $\Emm$.
It is clear that $\ciagi{\embes\Emm}$ is equivalent to $\ovr\Emm$ with embeddings and $\ciagi{(\embes\Emm,\homos\Emm)}$ is equivalent to $\ovr\Emm$ with all homomorphisms.

Recall that a model $X\in\ovr\Emm$ is {\em algebraically closed} if for every formula
$$\phi(x_0,\dots,x_{k-1}, y_0,\dots,y_{\ell-1})$$
that is a finite conjunction of atomic formulae, for every $a_0,\dots, a_{k-1} \in X$, if there exists an extension $X'\sups X$ in $\ovr\Emm$ satisfying
$$X'\models (\exists\;y_0,\dots,y_{\ell-1})\; \phi(a_0,\dots,a_{k-1},y_0,\dots,y_{\ell-1})$$
then there exist $b_0, \dots, b_{\ell-1} \in X$ such that
$X\models \phi(a_0,\dots,a_{k-1},b_0,\dots,b_{\ell-1})$.

\begin{prop}\label{ppanfgaju}
Let $\Emm$ be a class of finitely generated models of a fixed first-order language.
Every $\Emm$-injective model in $\ovr\Emm$ is algebraically closed.
\end{prop}

\begin{pf}
Fix an $\Emm$-injective model $X \in \ovr\Emm$.
Fix $X'\sups X$ and assume
$X' \models (\exists\; \vec y)\; \phi(\vec a, \vec y)$ for some $k$-tuple $\vec a$ of elements of $X$, where $\phi(\vec x,\vec y)$ is a finite conjunction of atomic formulae and $\vec x, \vec y$ are shortcuts for $(x_0,\dots, x_{k-1})$ and $(y_0,\dots,y_{\ell-1})$, respectively.

Let $S\in \Emm$ be a submodel of $X$ that contains $\vec a$.
Let $T\in \Emm$ be a submodel of $X'$ containing $S$ and a fixed tuple $\vec b$ such that $X' \models \phi(\vec a, \vec b)$.
Then also $T \models \phi(\vec a, \vec b)$, because this property is absolute for $\phi$.
Using the $\Emm$-injectivity of $X$, find a homomorphism $\map fTX$ satisfying $f\rest S = \id S$.
Finally, let $\vec c = (f(b_0),\dots,f(b_{\ell-1}))$, where $\vec b = (b_0,\dots,b_{\ell-1})$.
Since $f$ is a homomorphism and $\phi$ is a conjunction of atomic formulae, we have that $X \models \phi(\vec a, \vec c)$.
\end{pf}

We shall say that a structure $M$ is \emph{$n$-generated} if there exists $S\subs M$ such that $|S|\loe n$ and $S$ generates $M$, that is, no proper submodel of $M$ contains $S$.
Recall that a first-order language is \emph{finite} if it contains finitely many predicates (constant, relation and function symbols).

\begin{prop}\label{pnfgju}
Let $\Emm$ be a class of finite models of a fixed first-order language $L$. 
Assume that either $L$ is finite or for every $\ntr$ there exist finitely many isomorphic types of $n$-generated models in $\Emm$.
Assume furthermore that $\Emm$ has the mixed amalgamation property. Then every algebraically closed $L$-model $X \in \ovr\Emm$ is $\Emm$-injective.
\end{prop}

\begin{pf}
Fix $S,T \in \Emm$ such that $S$ is a submodel of $T$.
Fix a homomorphism $\map fSX$.
Using the mixed amalgamation, we can find an extension $X'\in \ovr\Emm$ of $X$ and a homomorphism $\map {f'}T{X'}$ such that $f'\rest S = f$.
Let $\Gee$ be the set of all functions $\map gTX$ satisfying $g \rest S = f$.
We need to show that some $g\in\Gee$ is a homomorphism.

Suppose first that there exist only finitely many $|T|$-generated structures in $\Emm$ and let $\En \subs \Emm$ be a finite set that contains isomorphic types of all of them.

Given $g\in \Gee$, denote by $g'$ a fixed isomorphism from the submodel generated by $\img gT$ onto a fixed model from the collection $\En$.
Note that $g$ is a homomorphism if and only if $g' \cmp g$ is a homomorphism.
Now observe that the set $\Ha = \setof{g' \cmp g}{g\in \Gee}$ is finite.

Let $S = \sett{s_i}{i < k}$ and $T\setminus S = \sett{t_j}{j < \ell}$.
Fix $g\in \Gee$ and suppose it is not a homomorphism.
There exists either a relation $R$ or a function $F$ and a finite sequence of elements of $T$ that witness this fact.
Let $\psi_g$ be an atomic formula describing this fact.
We may assume that $\psi_g$ has $k+\ell$ free variables, the first $k$ are supposed to denote $s_0,\dots,s_{k-1}$ and the latter ones 
$t_0,\dots,t_{\ell-1}$.
Let $\phi$ be the conjunction of all formulae $\psi_g$, where $g\in\Gee$.
Then $T\models \phi(s_0, \dots, s_{k-1}, t_0, \dots, t_{\ell-1})$ and, since $f'$ is a homomorphism,
\begin{equation}
X' \models \phi \Bigl(f(s_0), \dots, f(s_{k-1}),  f'(t_0), \dots, f'(t_{\ell-1})\Bigr).
\tag{1}\label{eqmfutione}
\end{equation}
Using the fact that $X$ is algebraically closed, find $\vec u = (u_0,\dots,u_{\ell-1})$ in $X$ such that
\begin{equation}
X \models \phi \Bigl(f(s_0), \dots, f(s_{k-1}), u_0, \dots, u_{\ell-1} \Bigr).
\tag{2}\label{eqmfugotwo}
\end{equation}
Let $g\in \Gee$ be such that $g(t_j) = u_j$ for $j < \ell$.
Then $g$ is a homomorphism. Indeed, otherwise there would be a witness (a relation or a function, plus some elements of $T$) saying that $g' \cmp g$ is not a homomorphism; however $\phi$ ``knows" all these witnesses, which gives rise to a contradiction.

Suppose now that $L$ is finite and consider again the set $\Gee$.
For each $g\in \Gee$, if $g$ is not a homomorphism, this is witnessed by an atomic formula $\psi_g$ and some elements of $T$. Now, even though the set $\Gee$ may be infinite, the number of atomic formulae with parameters in $T$ is finite.
As before, let $\phi(\vec x, \vec y)$ collect all of them.
Again, $X \models \phi(\vec s, \vec t)$ and consequently (\ref{eqmfutione}) holds.
Since $X$ is algebraically closed, we can find $\vec u$ such that (\ref{eqmfugotwo}) holds.
Finally, $g\in\Gee$ satisfying $g(t_j) = u_j$ ($j < \ell$) is the desired homomorphism.
\end{pf}

Following Dolinka~\cite{Dolinka}, we say that a class of models $\Emm$ has the \emph{1-point homomorphism extension property} (briefly: 1PHEP) if for every embedding $\map iAB$ and for every surjective homomorphism $\map fAC$, where $A,B,C \in \Emm$ and $B$ is generated by $A\cup \sn b$ for some $b\in B$, there exist an embedding $\map jCD$ and a homomorphism $\map gBD$ for which the diagram
$$\xymatrix{
C \ar@{ >->}[r]^j & D \\
A \ar@{->>}[u]^f \ar@{ >->}[r]_i & B \ar@{->>}[u]_g
}$$
commutes.
Let us say that an embedding $\map iAB$ is \emph{primitive} if $B$ is generated by one element from $\img iA$.
Clearly, every embedding is a composition of primitive embeddings.
Furthermore, every homomorphism is the composition of a surjective homomorphism and an embedding.
These facts, together with easy induction (see Proposition~\ref{pprimitivv}) show that 1PHEP is equivalent to the mixed amalgamation property of $\pair{\embes \Emm}{\homos \Emm}$.

Combining Theorem~\ref{tmejnn}, Propositions~\ref{ppanfgaju},~\ref{pnfgju} and the remarks above, we obtain a strengthening of Dolinka's result \cite{Dolinka}:

\begin{wn}
Let $\Emm$ be a \fra\ class of finite models of a fixed first-order language $L$.
Assume that $L$ is finite or for every $\ntr$ the number of isomorphism types of $n$-generated structures in $\Emm$ is finite.
Assume further that $\Emm$ has the pushout property and the 1PHEP.
Let $U\in\ovr\Emm$ be the \fra\ limit of $\Emm$.
For a model $X\in \ovr\Emm$ the following conditions are equivalent.
\begin{enumerate}
	\item[(a)] $X$ is a retract of $U$.
	\item[(b)] $X$ is algebraically closed.
\end{enumerate}
\end{wn}

The ``pushout property" in the statement above means that $\embes\Emm$ has pushouts in $\homos\Emm$. This assumption may of course be replaced by a weaker one, namely, that $\pair{\embes\Emm}{\homos\Emm}$ has the amalgamated extension property.

Note that a \fra\ class of finite models of a finite language may fail the condition concerning the number of $n$-generated models.
For example, let $L$ consist of a unique unary function symbol $P$ and let $\Emm$ be the class of all finite $L$-models.
That is, every model $S \in \Emm$ is endowed with a function $\map {P^S}SS$ and $\map fST$ is a homomoprhism iff $f(P^S(x)) = P^T(f(x))$ for every $x\in S$.
It is an easy exercise to check that $\Emm$ is a \fra\ class with the mixed pushout property, therefore the corollary above applies.
On the other hand, for each $\ntr$ there exists a $1$-generated structure $S_n \in \Emm$ of cardinality $n$. Namely, $S_n = \{0,\dots,n-1\}$ with the function $P$ defined by $P(n-1)=0$ and $P(i) = i+1$ for $i < n-1$.
Thus, there are infinitely many $1$-generated structures in $\Emm$.
Note that a countable $L$-structure $\pair XP$ belongs to $\ovr \Emm$ if and only if for every finite set $A\subs X$ there exists a finite set $S\subs X$ such that $A\subs S$ and $\img PS \subs S$.

There are some natural \fra\ classes of finite models of infinite languages and with infinitely many $n$-generated structures, for which Proposition~\ref{pnfgju} (and consequently the corollary above) still hold.
In section~\ref{sstegbjer} below, we shall investigate \fra\ classes of metric spaces, showing that the possibility of characterizing injectivity by ``being algebraic closed" depends on the language specifying the objects.

\subsection{A note on homomorphism-homogeneous structures}

In connection with (classical model-theoretic) \fra\ limits, there is an interesting notion of homomorphism-homogeneous structures, introduced recently by Cameron and Ne\v set\v ril~\cite{CamNes} and studied already by several authors (see \cite{Mas07}, \cite{Ilic}, \cite{CamLock}, \cite{RusSch}, \cite{MasPech}).
Namely, a (usually countable) structure $M$ is \emph{homomorphism-homogeneous} if every homomorphism between its finitely-generated substructures extends to an endomorphism of $M$.
It is clear that this notion can be defined in category-theoretic language, using a pair of categories $\pair \fK\fL$ as before, where $\fL$-arrows mean ``homomorphisms" and $\fK$-arrows mean ``embeddings".
It turns out that homomorphism-homogeneity is strictly related to injectivity, as we show below.

\begin{df}
Fix two categories $\fK\subs \fL$ with the same objects.
We say that an object $X \in \ob{\ciagi\fK}$ is \emph{$\fL$-homogeneous in $\ciagi{(\fK,\fL)}$} if for every $\ciagi\fK$-arrow $\map jaX$ such that $a\in \ob\fK$, for every $\ciagi{(\fK,\fL)}$-arrow $\map faX$, there exists a $\ciagi{(\fK,\fL)}$-arrow $\map FXX$ satisfying $F \cmp j = f$. This is described in the diagram below.
$$\xymatrix{
a \ar@{ >->}[rr]^j \ar[rrd]_f & & X \ar@{-->}[d]^F \\
& & X
}$$
\end{df}

Note that the arrow $f$ is of the form $x_n^\infty \cmp f'$ for some $f' \in \fL$.
That is why the definition above really speaks about $\fL$-homogeneity, not $\ciagi{(\fK,\fL)}$-homogeneity.
It can actually be viewed as a variation on the mixed amalgamation property, which is witnessed by the results below.

\begin{lm}\label{lnesetrilrtn}
Let $\pair \fK\fL$ be a pair of categories such that $\fK \subs \fL$ and let $X,Y \in \ob{\ciagi\fK}$ be such that $X$ is $\fK$-injective in $\ciagi{(\fK,\fL)}$. Then for every $\ciagi\fK$-arrow $\map jaY$ with $a\in \ob\fK$, for every $\ciagi{(\fK,\fL)}$-arrow $\map faX$, there exists a $\ciagi{(\fK,\fL)}$-arrow $\map FYX$ for which the diagram
$$\xymatrix{
a \ar@{ >->}[rr]^j \ar[rrd]_f & & Y \ar[d]^F \\
& & X
}$$
commutes.
\end{lm}

\begin{pf}
The arrow $j$ factorizes through some $y_k$, that is, $j = y_k^\infty \cmp i$ for some $\fK$-arrow $\map ia{y_k}$.
Using $\fK$-injectivity, we construct inductively $\fL$-arrows $\map {f_n}{y_n}X$ for $n \goe k$ so that $f_k \cmp i = f$ and $f_{n+1} \cmp y_n^{n+1} = f_n$ for $n > k$.
This gives rise to an arrow of sequences $F = \sett{f_n}{n \goe k}$ satisfying $F \cmp j = f$.
\end{pf}

Letting $X = Y$ in the lemma above, we obtain:

\begin{wn}\label{wnrgtboh}
Let $\fK \subs \fL$ be a pair of categories. Every $\fK$-injective object is $\fL$-homogeneous in $\ciagi{(\fK,\fL)}$.
\end{wn}

The equivalence (b)$\iff$(c) in the next statement, in the context of model theory, has been noticed by Dolinka~\cite[Prop. 3.8]{DolinkaBerg}.

\begin{prop}
Let $\fK \subs \fL$ be a pair of categories and let $\fK$ have a \fra\ sequence $U \in \ob{\ciagi\fK}$.
The following properties are equivalent:
\begin{enumerate}
	\item[(a)] $U$ is $\fK$-injective in $\ciagi{(\fK,\fL)}$.
	\item[(b)] $U$ is $\fL$-homogeneous in $\ciagi{(\fK,\fL)}$.
	\item[(c)] $\pair \fK\fL$ has the mixed amalgamation property.
\end{enumerate}
\end{prop}

\begin{pf}
Implication (c)$\implies$(a) has been proved in  Proposition~\ref{pwkindrzektif}.
Implication (a)$\implies$(b) is a consequence of Corollary~\ref{wnrgtboh}.
It remains to show that (b)$\implies$(c).

Suppose $U$ is $\fL$-homogeneous in $\ciagi{(\fK,\fL)}$ and fix a $\fK$-arrow $\map jca$ and an $\fL$-arrow $\map fcb$.
Using the property of being a \fra\ sequence, find $\fK$-arrows $\map ia{u_k}$ and $\map eb{u_\ell}$ with some $k,\ell < \nat$.
Since $U$ is $\fL$-homogeneous, there exists a $\ciagi{(\fK,\fL)}$-arrow $\map FUU$ satisfying $F \cmp u_k^\infty \cmp i \cmp j = u_\ell^\infty \cmp e \cmp f$.
Finally, find an $\fL$-arrow $\map g{u_k}{u_m}$ with $m > \ell$, such that $u_m^\infty \cmp g = F \cmp u_k^\infty$.
The situation is described in the following diagram.
$$\xymatrix{
c \ar@{ >->}[r]^j \ar[d]_f & a \ar@{ >->}[r]^i & u_k \ar[rd]^g \ar@{ >->}[rrr] & & & U \ar[d]^F \\
b \ar@{ >->}[rr]_e & & u_\ell \ar@{ >->}[r]_{u_\ell^m} & u_m \ar@{ >->}[rr] & & U
}$$
Thus, $j$ and $f$ are amalgamated by a $\fK$-arrow $u_\ell^m \cmp e$ and an $\fL$-arrow $g \cmp i$.
\end{pf}

Under certain natural assumptions, we are able to characterize homo\-mor\-phism-homo\-gen\-eous objects.
In the next statement we deal with countable categories, but what we really have in mind is the existence of countably many isomorphic types of arrows. 
For example, the category of finite sets is a proper class, yet it is obviously equivalent to a countable category.

\begin{tw}\label{tmejndwa}
Let $\fK \subs \fL$ be a pair of categories such that $\pair \fK\fL$ has the mixed pushout property, $\fL$ is countable, and $\fK$ has the initial object $0$.
For a sequence $X \in \ob{\ciagi\fK}$, the following properties are equivalent.
\begin{enumerate}
	\item[(a)] $X$ is $\fL$-homogeneous in $\ciagi{(\fK,\fL)}$.
	\item[(b)] There exists a subcategory $\fK_0$ of $\fK$ such that $0$ is initial in $\fK_0$, $X \in \ob{\ciagi{\fK_0}}$, $\pair{\fK_0}\fL$ has the mixed pushout property, and $X$ is $\fK_0$-injective in $\ciagi{(\fK_0,\fL)}$.
	\item[(c)] There exists a subcategory $\fK_0$ of $\fK$ such that $0$ is initial in $\fK_0$, $X \in \ob{\ciagi{\fK_0}}$, $\pair{\fK_0}\fL$ has the mixed pushout property, and $X$ is a retract of a \fra\ sequence in $\fK_0$.
\end{enumerate}
\end{tw}

The existence of the initial object in $\fK$ is not essential, but to remove it we would have to make more technical assumptions involving the joint embedding property.

\begin{pf}
The equivalence (b)$\iff$(c) is contained in Theorem~\ref{tmejnn}.
The fact that $\fK_0$ is countable has been used here for the existence of a \fra\ sequence.
Implication (b)$\implies$(a) is contained in Corollary~\ref{wnrgtboh}.
It remains to show that (a)$\implies$(b).

We may assume that $x_0 = 0$ in the sequence $X$.
Let $\Es = \setof{x_n^m}{n \loe m, \; n,m\in\nat}$.
Then $\Es$ is a subcategory of $\fK$ that contains the intial object $0$.
We first check that $X$ is $\Es$-injective.
Fix a $\ciagi{(\fK,\fL)}$-arrow $\map f{x_n}X$ and fix $m > n$.
Since $X$ is $\fL$-homogeneous, there exists a $\ciagi{(\fK,\fL)}$-arrow $\map FXX$ satisfying $F \cmp x_n^\infty = f$.
Note that $x_n^\infty = x_m^\infty \cmp x_n^m$, therefore $(F \cmp x_m^\infty) \cmp x_n^m = f$, which shows the $\Es$-injectivity of $X$.

Now let $\fK_0$ consist of all $\fK$-arrows $\map jca$ such that $X$ is $j$-injective in $\ciagi{(\fK,\fL)}$ and there exists at least one $\ciagi{(\fK,\fL)}$-arrow from $c$ to $X$.
That is, for every $\ciagi{(\fK,\fL)}$-arrow $\map fcX$, there exists a $\ciagi{(\fK,\fL)}$-arrow $\map gaX$ satisfying $g \cmp j = f$.
The second assumption is needed for keeping $0$ initial in $\fK_0$, namely, $X$ should also be injective for the (unique) arrow $0\to a$.
It is clear that $\fK_0$ is a subcategory of $\fK$ containing $\Es$.
In particular, $X\in \ob{\ciagi {\fK_0}}$.
It remains to show that $\pair{\fK_0}\fL$ has the mixed pushout property.
For this aim, fix an $\fK_0$-arrow $\map jcb$, an $\fL$-arrow  $\map pca$, and let $\map kaw$, $\map qbw$ be such that $k\in\fK$, $q\in\fL$ and
$$\xymatrix{
b \ar[r]^q & w \\
c \ar@{ >->}[u]^j \ar[r]_p & a \ar@{ >->}[u]_k
}$$
is a pushout square in $\fL$.
Fix a $\ciagi{(\fK,\fL)}$-arrow $\map faX$.
Since $X$ is $j$-injective, there exists a $\ciagi{(\fK,\fL)}$-arrow $\map gbX$ satisfying $g \cmp j = f \cmp p$.
Both arrows $f$ and $g$ are factorized through some $x_n$, namely,
$f = x_n^\infty \cmp f'$ and $g = x_n^\infty \cmp g'$ for some $\fL$-arrows $f', g'$.
Using the property of a pushout, we find a unique $\fL$-arrow $\map hw{x_n}$ satisfying $h \cmp k = f'$ and $h \cmp q = g'$.
In particular, $\ovr f = x_n^\infty \cmp h$ has the property that $\ovr f \cmp k = f$.
This shows that $X$ is $k$-injective in $\ciagi{(\fK,\fL)}$ and completes the proof.
\end{pf}

Unfortunately, the result above is not fully applicable to \fra\ classes. Namely, in case where $\fK$ is a countable \fra\ class, $\fK_0$ may not be a full subcategory of $\fK$.
This is demonstrated below, for the class of finite graphs.

\begin{ex}
Let $X$ be the two-element complete graph.
It is clear that $X$ is homomorphism-homogeneous (and also ultrahomogeneous).
We consider graphs without loops, therefore every endomorphism of $X$ is an automorphism.
More precisely, we consider the pair $\pair\fK\fL$, where $\ob\fK = \ob\fL$ are all finite simple graphs, the $\fK$-arrows are embeddings and the $\fL$-arrows are graph homomorphisms.

Let $\fK_0$ be any subcategory of $\fK$ that has pushouts in $\fL$ and contains all the embeddings of subgraphs of $X$.
So, $\ob{\fK_0}$ contains the empty graph and complete subgraphs of size $\loe2$. 
The pushout with embeddings of the empty graph is just the coproduct (disjoint sum), there $\ob{\fK_0}$ contains the 2-element graph $D$ with no edges.
Furthermore, $\ob{\fK_0}$ contains the graph $G$ whose set of vertices is $\{-1,0,1\}$ and the edges are $\dn{-1}0$ and $\dn01$.
Such a graph comes from the pushout of two embeddings of the one-element graph into $X$.
Now consider an embedding $\map jDG$ such that $\img jD = \dn{-1}1$.
Let $\map fDX$ be one-to-one.
Clearly, $f$ is a homomorphism and no homomorphism $\map gGX$ satisfies $g\cmp j  = f$.
This shows that $X$ is not $j$-injective. In particular, $\fK_0$ is not a full subcategory of $\fK$.
\end{ex}

\subsection{Metric spaces}\label{sstegbjer}

We shall now discuss a concrete model-theoretic application of our result: Retracts of the universal metric space of Urysohn.
Let $\metrics$ be the category of finite metric spaces with isometric embeddings.
The objects of $\metrics$ are models of a first-order language:
For each $r>0$ we can define the binary relation $D_r(x,y) \iff d(x,y) < r$, where $d$ denotes the metric on a fixed set $X$. The axioms of a metric can be rephrased in terms of the relations $D_r$.
For example, the triangle inequality follows from the following (infinitely many) formulae:
$$D_r(x,z) \land D_s(z,y) \implies D_{r+s}(x,y).$$
Note that it suffices to consider the relations $D_r$ with $r$ positive rational: The metric is then defined by $d(x,y) = \inf_{r \in \Qyu^+}{D_r(x,y)}$, where $\Qyu^+$ denotes the set of all positive rationals.
In other words, metric spaces can be described in a countable language.
It is clear that, in this language, a homomorphism of metric spaces is a non-expansive map. Recall that $\map fXY$ is \emph{non-expansive} if $d_Y(f(p),f(q)) \loe d_X(p,q)$ for every $p,q \in X$, where $d_X$, $d_Y$ denote the metrics on $X$ and $Y$ respectively.

It is also possible to describe a metric space by similar relations $D_r$, now meaning that the distance is $\loe r$. We shall see later that, even though both languages describe the same objects, the notion of being algebraically closed is completely different.

Clearly, the language of metric spaces is infinite and there exist infinitely many types of $2$-element metric spaces (even when restricting to rational distances), therefore one cannot apply Dolinka's result here. Moreover, $\metrics$ is formally not a \fra\ class, because it contains continuum many pairwise non-isomorphic objects. It becomes a \fra\ class when restricting to spaces with rational distances. However, in that case we cannot speak about complete metric spaces.
In any case, our main result is applicable to the complete metric space of Urysohn, as we show below.

The following lemma, in a slightly different form, can be found in \cite[Lemma 3.5]{DolMas}.

\begin{lm}\label{lnjgetrgg}
Let $\map fXY$ be a non-expansive map of nonempty finite metric spaces. Assume $X\cup\sn a$ is a metric extension of $X$. Then there exists a metric extension $Y\cup \sn b$ of $Y$ such that
$$\xymatrix{
Y \ar@{ >->}[r] & Y\cup\sn b \\
X \ar[u]^f \ar@{ >->}[r] & X\cup\sn a \ar[u]_g
}$$
where $g\rest X = f$ and $g(a) = b$, is a pushout square in the category of metric spaces with non-expansive maps.
Furthermore
\begin{equation}
d(y,b) = \min_{x\in X} \Bigl( d(y,f(x)) + d(x,a) \Bigr)
\tag{M}\label{eqemfwoef}
\end{equation}
for every $y\in Y$.
\end{lm}

The statement obviously fails when $X = \emptyset$ and $Y \nnempty$.

\begin{pf}
We first need to show that (\ref{eqemfwoef}) defines a metric on $Y\cup \sn b$.
Of course, only the triangle inequality requires an argument.
Fix $y,y_1 \in Y$.
Find $x_1\in X$ such that $d(y_1,b) = d(y_1,f(x_1)) + d(x_1,a)$.
Using the triangle inequality in $Y$, we get
$$d(y,b) \loe d(y,f(x_1)) + d(x_1,a) \loe d(y,y_1) + d(y_1,f(x_1)) + d(x_1,a) = d(y,y_1) + d(y_1,b).$$
Now find $x\in X$ such that $d(y,b) = d(y,f(x)) + d(x,a)$.
Using the triangle inequalities in $X$ and $Y$, and the fact that $d(f(x),f(x_1)) \loe d(x,x_1)$, we obtain
\begin{align*}
d(y,b) + d(y_1,b) &= d(y,f(x)) + d(x,a) + d(y_1,f(x_1)) + d(x_1,a) \\
&\goe d(y,f(x)) + d(x,x_1) + d(y_1,f(x_1)) \\
&\goe d(y,f(x)) + d(f(x),f(x_1)) + d(y_1,f(x_1)) \\
&\goe d(y,y_1).
\end{align*}
Thus, $d$ defined by (\ref{eqemfwoef}) fulfills the triangle inequality.

Given $x\in X$, we have $d(g(x),g(a)) = d(f(x),b) \loe d(f(x),f(x)) + d(x,a) = d(x,a)$. This shows that $g$ is non-expansive.

Finally, assume $\map p{X\cup\sn a}W$ and $\map qYW$ are non-expansive maps such that $p\rest X = q \cmp f$.
We need to show that there exists a unique non-expansive map $\map h{Y\cup\sn b}W$ satisfying $h \cmp g = p$ and $h \rest Y = q$.
The uniqueness of $h$ is clear, namely $h(b) = h(g(a)) = p(a)$.
It remains to verify that $h$ is non-expansive.

Suppose otherwise and fix $y\in Y$ such that $d(h(y),h(b)) > d(y,b)$.
Find $x\in X$ such that $d(y,b) = d(y,f(x)) + d(x,a)$.
So we have
\begin{equation}
d(h(y),p(a)) > d(y,f(x)) + d(x,a).
\tag{*}\label{eqgwiazz}
\end{equation}
Knowing that $p$ and $q$ are non-expansive, we get
\begin{equation}
d(p(x), p(a)) \loe d(x,a) \oraz d(q(y),q(f(x)) \loe d(y,f(x)).
\tag{**}\label{eqgwizdd}
\end{equation}
Note that $q(f(x)) = p(x)$ and $q(y) = h(y)$. Finally, (\ref{eqgwiazz}) and (\ref{eqgwizdd}) give
$$d(h(y),p(a)) > d(p(x),p(a)) + d(h(y),p(x))$$
which contradicts the triangle inequality in $W$.
This completes the proof.
\end{pf}

We say that a metric space $\pair Xd$ is \emph{finitely hyperconvex} if for every finite family of closed balls
$$\Aaa = \left\{ \clbal(x_0,r_0), \clbal(x_1,r_1), \dots, \clbal(x_{n-1},r_{n-1}) \right\}$$
such that $\bigcap \Aaa = \emptyset$, there exist $i,j < n$ such that
$$d(x_i,x_j) > r_i + r_j.$$
This is a weakening of the notion of a \emph{hyperconvex metric space}, due to Aronszajn \& Panitchpakdi~\cite{AroPan}, where the family above may be of arbitrary cardinality.
Actually, the authors of \cite{AroPan} had already considered $\kappa$-hyperconvex metric spaces; finite hyperconvexity corresponds to $\aleph_0$-hyperconvexity.
A variant of finite hyperconvexity (with closed balls replaced by open balls) has been recently studied by Niemiec~\cite{NiemiecANR} in the context of topological absolute retracts.

The following facts relate this definition to our main topic.
The first one should be well known to readers familiar with hyperconvexity, namely, every metric space embeds isometrically into a hyperconvex one.

\begin{lm}\label{lhtrgoo}
Let $X$ be a finite metric space and let $\Aaa = \sett{\clbal(x_i,r_i)}{i<N}$ be a family of closed balls such that $N\in\nat$ and $d(x_i,x_j) \loe r_i + r_j$ for every $i,j<N$. Then there exists a metric extension $X\cup \sn a$ of $X$ such that $d(a,x_i) \loe r_i$ for every $i<N$.
\end{lm}

\begin{pf}
Fix $a\notin X$ and define
\begin{equation}
d(a,x) = \min_{i < N}\Bigl( d(x,x_i) + r_i \Bigr).
\tag{*}\label{eqnstearr}
\end{equation}
Obviously, $d(a,x_i)\loe r_i$. It remains to check that (\ref{eqnstearr}) indeed defines a metric on $X\cup\sn a$.
It is the triangle inequality that requires a proof.
Fix $x,y\in X$ and fix $k < N$ such that $d(a,y) = d(y,x_k) + r_k$.
Then
$$d(a,x) \loe d(x,x_k) + r_k \loe d(x,y) + d(y,x_k) + r_k = d(x,y) + d(y,a).$$
This shows that $d(x,a) \loe d(x,y) + d(y,a)$.
Now fix $i < N$ such that $d(a,x) = d(x,x_i) + r_i$. We have that $d(x_i,x_k) \loe r_i + r_k$, therefore
\begin{align*}
d(a,x) + d(a,y) &= d(x,x_i) + r_i + d(y,x_k) + r_k \\
&\goe d(x,x_i) + d(y,x_k) + d(x_i,x_k) \goe d(x,y).
\end{align*}
This shows that $d$ defined by (\ref{eqnstearr}) satisfies the triangle inequality.
\end{pf}

The next lemma is a special case of two results of Aronszajn \& Panitchpakdi, namely, Theorem 2 on page 413 and Theorem 3 on page 415 in \cite{AroPan}.
We present the proof for the sake of completeness.

\begin{lm}\label{ltebhr}
A metric space is finitely hyperconvex if and only if it is injective with respect to isometric embeddings of finite metric spaces.
\end{lm}

\begin{pf}
Let $X$ be a finitely hyperconvex metric space and fix a non-expansive map $\map fSX$, where $S$ is a finite metric space.
It suffices to show that $f$ can be extended to a non-expansive map $\map {f'}TX$ whenever $T$ is a metric extension of $S$ and $T\setminus S = \sn a$.
Fix $T = S \cup \sn a$ and let
$$\Aaa = \setof {\clbal(f(s), r_s)}{s\in S},$$
where $r_s = d(s,a)$.
Given $s,s_1\in S$, we have that $d(f(s),f(s_1)) \loe d(s,s_1) \loe r_s + r_{s_1}$.
Since $X$ is finitely hyperconvex, there exists $b\in \bigcap\Aaa$.
This means that $d(b,f(s)) \loe d(s,a)$ for every $s\in S$.
Thus, setting $f'(a) = b$ and $f'\rest S = f$, we obtain a non-expansive extension of $f$.
This shows the ``only if" part.

For the ``if" part, fix a family $\Aaa = \sett{\clbal(x_i, r_i)}{i < N}$ in $X$, so that $d(x_i, x_j) \loe r_i + r_j$ for $i,j < N$.
Let $S = \{ x_0,x_1, \dots, x_{N-1} \}$ and endow $S$ with the metric inherited from $X$.
Let $T = S\cup \sn a$ be a metric extension of $S$ such that $d(a,x_i) \loe r_i$ for $i < N$.
It exists by Lemma~\ref{lhtrgoo}.
Applying the injectivity of $X$, we can find a non-expansive extension $\map gTX$ of the inclusion $S\subs X$.
Let $b = g(a)$. Then $d(b,x_i) \loe d(a,x_i) \loe r_i$ for $i < N$.
This shows that $\bigcap \Aaa \nnempty$.
\end{pf}

\begin{tw}\label{tfndgbojf}
Given a Polish space $\pair Xd$, the following properties are equivalent:
\begin{enumerate}
	\item[(a)] $\pair Xd$ is a non-expansive retract of the universal Urysohn space $\U$.
	\item[(b)] $\pair Xd$ is finitely hyperconvex.
	\item[(b')] $\pair Xd$ is injective with respect to isometric embeddings of finite metric spaces.
\end{enumerate}
\end{tw}

\begin{pf}
The equivalence (b)$\iff$(b') is contained in Lemma~\ref{ltebhr}.

(a)$\implies$(b') 
Assume $X\subs \U$ and $\map r\U X$ is a non-expansive retraction.
Fix finite metric spaces $S \subs T$ and a non-expansive map $\map fSX$.
Using the mixed pushout property (a consequence of Lemma~\ref{lnjgetrgg} and Proposition~\ref{pjetghsd}), we can find an isometric embedding $\map j{\img fS}W$ and a non-expansive map $\map gTW$ such that $W$ is a finite metric space and $g\rest S = j \cmp f$.
Using the ultrahomogeneity of $\U$, we can find an isometric embedding $\map hW\U$ such that $h \cmp j$ is the inclusion $\img fS \subs \U$.
Finally, let $p = r \cmp h \cmp g$. Then $\map pTX$ is a non-expansive map and $p\rest S = f$.

(b')$\implies$(a)
Fix a Polish space $X$ satisfying (b').
Fix a countable dense set $D\subs X$.
Let $K_0 = \Qyu \cup \setof{d(x,y)}{x,y\in D}$ and let $K$ be the subsemigroup of $\pair\Err +$ generated by $K_0$.
Consider the category of nonempty finite metric spaces with distances in $K$ (we call them \emph{$K$-metric spaces}). This category is countable, therefore it has a \fra\ sequence. This \fra\ sequence defines a countable metric space $E$ whose completion is, by uniqueness, the Urysohn space.
Enlarging $D$ to a countable set, we may assume that it is injective with respect to isometric embeddings of finite $K$-metric spaces.
By Theorem~\ref{tmejnn}, $D$ is a non-expansive retract of $E$ and consequently $X$ is a non-expansive retract of $\U$.
\end{pf}

Let us note that in the statement above only the implication (b)$\implies$(a) appears to be new, the other ones are standard arguments easily adapted from \cite{AroPan}.
The main ingredient needed here is the fact that Urysohn's space is finitely hyperconvex, which follows directly from Lemma~\ref{ltebhr} above.

It is easy to see that the results above remain valid for the bounded version of the Urysohn space, called the \emph{Urysohn sphere}.
Denote by $\metrics_1$ the class of all finite metric spaces of diameter $\loe 1$.
There is an obvious functor mapping $\pair Xd \in \metrics$ to $\pair X{d_C} \in \metrics_1$, where
$$d_C(x,y) = \min\{d(x,y), C\}.$$
Applying this functor, we can easily conclude that Lemmata~\ref{lnjgetrgg}, \ref{lhtrgoo} and~\ref{ltebhr} hold for arbitrary classes of the form  $\metrics_1$.
However, we need to specify the more general version of hyperconvexity.
Namely, we say that $\pair Xd$ is \emph{finitely $1$-hyperconvex} if for every finite family of closed balls $\Bee = \sett{\clbal(x_i,r_i)}{i<n}$ with $r_i \loe 1$ for $i<n$, it holds that $\bigcap \Bee \nnempty$ whenever $d(x_i,x_j) \loe r_i + r_j$ for every $i,j < n$.

The bounded version of Theorem~\ref{tfndgbojf} is as follows.

\begin{tw}
Given a separable complete metric space $X$ of diameter $\loe 1$, the following conditions are equivalent.
\begin{enumerate}
	\item[(a)] $X$ is a non-expansive retract of the Urysohn sphere.
	\item[(b)] $X$ is finitely $1$-hyperconvex.
	\item[(c)] $X$ is injective with respect to isometric embeddings of finite metric spaces of diameter $\loe 1$.
\end{enumerate}
\end{tw}

The theorem above speaks about complete metric spaces, however
we can also formulate a version by restricting distances to a countable subsemigroup $S$ of $[0,+\infty)$.
In that case, the class $\Emm_S$ of finite metric spaces with distances in $S$ is countable and we can consider its \fra\ limit $U_S\in \ovr \Emm_S$, a (possibly non-complete) countable ultrahomogeneous $S$-metric space.
By the remarks above, we conclude that $X \in \ovr \Emm_S$ is a non-expansive retract of $U_S$ if and only if it is finitely $S$-hyperconvex (with the obvious meaning of $S$-hyperconvexity).
This gives rise to the announced example showing that ``being algebraically closed" for metric spaces may or may not be equivalent to injectivity.

\begin{ex}
Consider the class $\Emm_\Qyu$ of finite rational metric spaces.
Assume that the language consists of relations $D_r$ ($r\in\Qyu$), where $D_r(x,y)$ means ``$d(x,y) < r$".
Using finite conjunctions of atomic formulae, there is no way to say that $X$ is finitely $\Qyu$-hyperconvex.
Indeed, consider $\Qyu$ as a metric space with the usual distance and take $X = \Qyu \setminus \sn 1$.
Clearly, $\Qyu$ is finitely $\Qyu$-hyperconvex, hence algebraically closed (see Proposition~\ref{ppanfgaju}).
Thus, $X$ is algebraically closed too, because of the strict inequalities in the relations $D_r$.
On the other hand, $X$ is obviously not $\Emm_\Qyu$-injective: The inclusion $\dn02 \subs X$ has no non-expansive extension onto $\{0,1,2\}$.

Finally, consider the same language for $\Emm_\Qyu$, but with a different interpretation. Namely, let $D_r(x,y)$ mean ``$d(x,y) \loe r$" ($r\in \Qyu$).
Now it is clear that ``being algebraically closed" implies ``being finitely $\Qyu$-hyperconvex", because of a version of Lemma~\ref{lhtrgoo} for rational metric spaces.
Thus, the two properties are equivalent and now it is true that a countable rational metric space is $\Emm_\Qyu$-injective if and only if it is algebraically closed.
\end{ex}

\subsection{Banach spaces}

Let $\banach$ denote the category of finite-dimensional Banach spaces (over the field of real or complex numbers) with linear transformations of norm $\loe 1$.
Let $\banachi$ denote the category of finite-dimensional Banach spaces with linear isometric embeddings.

The following lemma is well known.
For the proof we refer to~\cite{ACCGMud}.

\begin{lm}
$\pair\banachi\banach$ has the mixed pushout property.
\end{lm}

A Banach space $X$ is \emph{1-complemented} in $Y$ if $X \subs Y$ and there exists a projection $\map P Y Y$ (i.e. a linear operator satisfying $P\cmp P = P$) of norm $1$ and $\img P Y = X$.
A Banach space $E$ is \emph{almost $1$-injective} for finite-dimensional spaces if, given finite-dimensional spaces $X \subs Y$, given a linear operator $\map T X E$ with $\norm T \loe 1$, given $\eps > 0$, there exists a linear operator $\map {\til T} Y E$ such that $\til T \rest X = T$ and $\norm{\til T} \loe 1 + \eps$.
The \emph{Gurarii space}~\cite{Gurarii} is a separable Banach space $\gur$ satisfying the following condition:
Given $\eps > 0$ and finite-dimensional spaces $X \subs Y$, every isometric embedding $\map e X \gur$ extends to an $\eps$-isometric embedding $\map {\til e} Y \gur$ (that is, $\til e$ is one-to-one and $\norm {\til e} \loe 1 + \eps$, $\norm {\til e^{-1}} \loe 1 + \eps$).
The fact that $\gur$ is unique up to a linear isometry was proved by Lusky~\cite{Lusky}; an elementary argument has been found recently, see~\cite{KubSol}.

We now would like to apply Theorem~\ref{tmejnn}. The obstacle is that the category $\banachi$ is too big, it does not have a \fra\ sequence.
On the other hand, given a countable $S\subs \banachi$ there exists a countable $\fK\subs \banachi$ such that $S\subs \fK$ and $\fK$ has pushouts in $\banach$.
The category $\fK$ has a \fra\ sequence. If $S$ is ``rich enough" then this \fra\ sequence induces the Gurarii space $\gur$. 
This way we obtain the following result, originally due to Wojtaszczyk~\cite{Wojtaszczyk}.

\begin{tw}
Let $E$ be a separable Banach space. The following properties are equivalent.
\begin{enumerate}
	\item[(a)] $E$ is linearly isometric to a $1$-complemented subspace of the Gurarii space.
	\item[(b)] $E$ is almost $1$-injective for finite-dimensional Banach spaces.
	\item[(c)] $E$ is an isometric $L^1$ predual.
\end{enumerate}
\end{tw}

\begin{pf}
(a)$\implies$(b) 
By the mixed pushout property, it is straightforward to see that the Gurarii space is almost $1$-injective.
Clearly, this property is preserved by $1$-complemented subspaces.

(b)$\implies$(c) This is part of the main result of Lindenstrauss \cite{Lin64}.
In fact, it is proved in \cite[Thm. 6.1]{Lin64} that (c) is equivalent to almost $1$-injectivity for Banach spaces of dimension $\loe 4$.

(c)$\implies$(a)
A result of Lazar \& Lindenstrauss~\cite{LazLin66} says that there exists a chain 
$E_0\subs E_1 \subs E_2 \subs \dots$ of finite-dimensional subspaces of $E$ whose union is dense in $E$ and each $E_n$ is isometric to some $\ell^\infty_{k(n)}$. In fact, due to Michael \& Pe{\l}czy\'nski~\cite{MicPel}, one may assume that $k(n) = n$ for $\ntr$, although this is not needed here.

Let $\fL$ be a countable subcategory of $\banach$ that contains all inclusions $E_n\subs E_{n+1}$ and a fixed chain defining the Gurarii space. Enlarging $\fL$ by adding countably many arrows, we may assume that it is closed under mixed pushouts, that is, the pair $\pair \fK \fL$ has the mixed pushout property, where $\fK = \fL \cap \banachi$.

Let $\gur$ denote the Gurarii space. By the assumptions on $\fL$, we have that both $\gur$ and $E$ are objects of $\ciagi\fK$.
Now observe that $E$ is $\fK$-injective in $\ciagi{(\fK,\fL)}$.
Indeed, if $\map fAE$ is an arrow in $\ciagi{(\fK,\fL)}$, where $A, B \in\ob\fK$ are such that $A\subs B$, then $f$ is an isometric embedding of $A$ into some $E_n$ (by the definition of arrows between sequences). 
It is easy and well known that every space isometric to $\ell^\infty_m$ is $1$-injective for all Banach spaces.
Thus, $f$ can be extended to a linear isometry $\map {\ovr f}B{E_n}$.

We actually need one more assumption on $\fK$: namely that $\ovr f\in \fK$ whenever $f\in \fK$. This can be achieved by a standard closing-off argument.

Finally, Theorem~\ref{tmejnn} implies that $E$ is isometric to a $1$-complemented subspace of $\gur$.
\end{pf}

A non-separable version of the above result is actually much simpler and comes exactly as a particular case of the uncountable version of Theorem~\ref{tmejnn}:

\begin{tw}\label{Tetnrw}
Assume the continuum hypothesis. Let $\Ve$ be the unique Banach space of density $\aleph_1$ that is of universal disposition for separable spaces. A Banach space of density $\loe \aleph_1$ is isometric to a $1$-complemented subspace of $\Ve$ if and only if it is $1$-separably injective.
\end{tw}

Some explanations are needed here.
Namely, a Banach space $V$ is \emph{of universal disposition} for separable spaces if for every separable Banach spaces $X \subs Y$, every isometric embedding of $X$ into $V$ extends to an isometric embedding of $Y$ into $V$.
Our result from~\cite{Kubfra} says that, under the continuum hypothesis, there exists a unique Banach space $\Ve$ of density $\aleph_1$ and of universal disposition for separable spaces.
Extensions of this result can be found in~\cite{ACCGMud}, where more general constructions of spaces of universal disposition are presented.
It is shown there that $2^{\aleph_0}$ is the minimal density of a Banach space of universal disposition for separable spaces.
Finally, assuming the continuum hypothesis, the space $\bV$ is the \fra\ limit of separable Banach spaces with linear isometric embeddings.
The notion of being ``$1$-separably injective" has obvious meaning; it has been recently studied in~\cite{ACCGMsi}.
In this context, Theorem~\ref{Tetnrw} complements the results of~\cite{ACCGMsi}.

\subsection{Linear orders}

Let $\kappa$ be an infinite cardinal and let $\los\kappa$ denote the class of all linearly ordered sets of cardinality $< \kappa$.
A homomorphism of linearly ordered sets will be called an \emph{increasing map}.
As mentioned before, $\flos$ gives a natural example of a pair $\pair {\embes {\flos}}{\homos {\flos}}$ failing the pushout property.
However, we have the following

\begin{prop}
For every infinite cardinal $\kappa$, the pair $\pair {\embes {\los\kappa}}{\homos {\los\kappa}}$ has property \hahah.
\end{prop}

\begin{pf}
Condition \haha1 follows from \haha3, because $\embes{\los\kappa}$ has an initial object (the empty set) and $\homos{\los\kappa}$ has a terminal object, the $1$-element linearly ordered set.
It remains to show \haha2 and \haha3.

Call an embedding $\map jAB$ \emph{primitive} if $|B\setminus \img jA| \loe1$.
It is clear that every increasing embedding is the colimit of a transfinite sequence of primitive embeddings.
We shall use an uncountable version of Proposition~\ref{pprimitivv}, which can be easily proved by transfinite induction, using the fact that the category $\embes{\los\kappa}$ is $\kappa$-continuous in $\homos{\los\kappa}$.

Denote by $\Pee$ the class of all primitive embeddings in $\embes{\los\kappa}$.
Let us prove first that $\pair \Pee{\homos{\los\kappa}}$ has the amalgamated extension property (condition \haha3).
Fix linearly ordered sets $C,A,B$ such that $A = C \cup \sn a$ and $B = C \cup \sn b$.
Fix increasing maps $\map fAL$ and $\map gBL$ such that $f\rest C = g\rest C$.
Formally, we have to assume that $a\ne b$.
Let $W = A \cup B$.
We let $a < b$ if $f(a) < g(b)$; we let $a > b$ otherwise.
It is clear, using the compatibility of $f$ and $g$, that this defines a linear order on $W$, extending the orders of $A$ and $B$. 
The unique map $\map hWL$ satisfying $h\rest A = f$ and $h\rest B = g$ is increasing.
This shows \haha3.

Now fix linearly ordered sets $C,A,B$ such that $A = C \cup \sn a$ with $a\notin C$, and fix an increasing map $\map fCB$.
Let
$$L = \bigcup_{c < a}(\leftarrow, f(c)] \oraz R = \bigcup_{c > a}[f(c), \rightarrow),$$
where $(\leftarrow, x] = \setof{p}{p \loe x}$ and $[x,\rightarrow) = \setof{p}{p \goe x}$.
Note that $B = L \cup R$ and $L\cap R$ is either empty or a singleton.
Let $W = B \cup \sn w$, where either $w \in L\cap R$ or $w\notin B$ in case where $L\cap R = \emptyset$.
In the latter case, define $x < w$ and $w < y$ for $x\in L$, $y\in R$.
Define $\map gAW$ by setting $g(a) = w$ and $g \rest C = f$.
Clearly, $g$ is increasing and the inclusion $B\subs W$ is primitive.
This shows \haha2 and completes the proof.
\end{pf}

Note that every increasing embedding of finite linear orders is left-invertible.
Thus, we immediately obtain the following result.

\begin{wn}
Every countable linear order is order-isomorphic to an increasing retract of the set of rational numbers.
\end{wn}

Of course, this result can be proved directly, realizing that $X \cdot \Qyu$ with the lexicographic ordering is isomorphic to $\Qyu$, whenever $X$ is a countable linear order.
Note that this completely answers Question~10.6 from~\cite{McPhee}.

Passing to the uncountable case, let us note that $\los{\omega_1}$ has the \fra\ limit if and only if the Continuum Hypothesis holds. 
Denote this \fra\ limit by $\Qyu_{\omega_1}$. 
It is easy to check that a linearly ordered set $X$ of cardinality $\omega_1$ is injective for countable linear orders (isomorphic to $\Qyu_{\omega_1}$) if and only if
for every countable sets $A,B \subs X$ such that $a < b$ for 
$a\in A$, $b\in B$, there exists $x\in X$ such that $a \loe x \loe b$ ($a < x < b$) whenever $a\in A$, $b\in B$ (one of the sets $A$, $B$ may be empty).
For example, the closed unit interval $[0,1]$ satisfies this condition, therefore it can be embedded as an increasing retract of $\Qyu_{\omega_1}$.

\separator

We finish with some remarks on reversed \fra\ sequences.
General theory of reversed \fra\ limits of finite models (of a first-order language)
was developed in~\cite{IrSol}.
The idea comes just by considering the opposite category.
More specifically, fix a class $\Emm$ of finite models and consider the pair $\pair{\quos\Emm}{\homos\Emm}$, where $\quos\Emm$ is the category whose objects are elements of $\Emm$ and arrows are quotient maps.
Now property \hahah\ is defined by reversing the arrows in all the diagrams.
For example, amalgamation is replaced by ``reversed amalgamation" and pushouts are replaced by pullbacks.
Sequences are now contravariant functors and it is natural to consider their limits endowed with the topology, inherited from the product of finite sets. It is not hard to see that precisely the continuous homomorphisms are induced by arrows between sequences.
It is worth noting that if $\Emm$ is closed under finite products and substructures then $\quos\Emm$ has pullbacks in $\homos\Emm$. The pullback of two quotient maps $\map f X Z$, $\map g Y Z$ is provided by the structure
$$w = \setof{\pair st \in X \times Y}{f(s) = g(t)}.$$
Coming back to finite linear orders, consider the pair $\pair{\quos\flos}{\homos\flos}$.
It is straightforward to see that $\quos\flos$ has no pullbacks in $\homos\flos$.
On the other hand, it is easy and standard to check that this pair has (the reversed variant of) property \hahah.
Note that every increasing quotient of finite linearly ordered sets is right-invertible. Thus, all sequences in $\quos\flos$ are ``finitely projective".
It is clear that the inverse \fra\ limit of $\flos$ is the Cantor set endowed with the standard linear order.
Thus, using Theorem~\ref{tmejnn} (or, more precisely, Corollary~\ref{wnrhojht}), we obtain the following well known fact which belongs to the folklore.

\begin{wn}\label{wpveryastonn}
Every compact metric totally disconnected linearly ordered space is a continuous increasing retract of the standard Cantor set.
\end{wn}

Again, it is not hard to prove this fact directly, by showing that a metric compact totally disconnected linearly ordered space $K$ can be isomorphically embedded into the Cantor set and constructing the retraction ``manually". Note that the reversed \fra\ theory would only say that $K$ is a continuous increasing quotient of the Cantor set, however not all continuous increasing quotient maps of the Cantor set are right-invertible.

\subsection*{Acknowledgments}
The author is indebted to the anonymous referee for several helpful remarks, in particular for pointing out the reference~\cite{McPhee}.

\end{document}